\documentclass{article}

\usepackage{arxiv}

\usepackage[utf8]{inputenc} 
\usepackage[T1]{fontenc}    
\usepackage{hyperref}       
\usepackage{url}            
\usepackage{booktabs}       
\usepackage{amsfonts}       
\usepackage{nicefrac}       
\usepackage{microtype}      
\usepackage{lipsum}
\usepackage{bm}

\newtheorem{theorem}{Theorem}[section]
\newtheorem{lemma}[theorem]{Lemma}
\newtheorem{proposition}[theorem]{Proposition}

\newtheorem{remark}[theorem]{Remark}
\newtheorem{assumption}[theorem]{Assumption}
\newenvironment{proof}{\textbf{Proof} }

\newcommand{\qed}{\nobreak \ifvmode \relax \else
      \ifdim\lastskip<1.5em \hskip-\lastskip
      \hskip1.5em plus0em minus0.5em \fi \nobreak
      \vrule height0.75em width0.75em depth0.25em\fi}

\usepackage{amsmath, amssymb, epsfig, graphicx, float, url}
\usepackage{subcaption}

\usepackage{algorithm, algorithmicx, algpseudocode}
\usepackage[colorinlistoftodos,bordercolor=orange,backgroundcolor=orange!20,line
color=orange,textsize=scriptsize]{todonotes}
 \newcommand{\tb}{\textcolor{black}}

\numberwithin{equation}{section}

\newcommand{\R}{\mathbb{R}}

\newcommand{\A}{\mathbb{A}}

\newcommand{\N}{\mathbb{N}}
\newcommand{\X}{\mathbb{S}}

\newcommand{\E}{\mathbb{E}}
\newcommand{\EE}{\mathcal{E}}
\newcommand{\PP}{\mathbb{P}}

\newcommand{\NN}{\mathcal{N}}
\newcommand{\M}{\mathcal{M}}
\newcommand{\W}{\mathcal{W}}

\newcommand{\I}{\mathcal{I}}

\title{Robust Markov Decision Process: \\
Beyond Rectangularity}

\author{
  Vineet Goyal\\
  IEOR Department\\
  Columbia University\\
  \texttt{vg2277@columbia.edu} \\
   \And
  Julien Grand-Cl{\'e}ment\\
  IS/OM Department\\
  HEC Paris\\
  \texttt{grand-clement@hec.fr} \\
}

\begin{document}
\maketitle

\begin{abstract}
We consider a robust approach to address uncertainty in model parameters in Markov Decision Processes (MDPs), {\color{black}which are} widely used to model dynamic optimization in many applications. Most prior works consider the case where the uncertainty on transitions related to different states is uncoupled and the adversary is allowed to select the worst possible realization for each state unrelated to others, potentially leading to highly conservative solutions. On the other hand, the case of general uncertainty sets is known to be intractable. We consider a factor model for probability transitions where the transition probability is a linear function of a {\em factor matrix} that is uncertain and belongs to a {\em factor matrix uncertainty set}. This is a fairly general model of uncertainty in probability transitions, allowing the decision maker to model dependence between probability transitions across different states and it is significantly less conservative than prior approaches. We show that under an underlying rectangularity assumption, we can efficiently compute an optimal robust policy under the factor matrix uncertainty model. Furthermore, we show that there is an optimal robust policy that is deterministic, which is of interest from an interpretability standpoint. We also introduce the robust counterpart of important structural results of classical MDPs, including the maximum principle and Blackwell optimality, and we provide a computational study to demonstrate the effectiveness of our approach in mitigating the conservativeness of robust policies.
\end{abstract}

\section{Introduction.}
\label{sec:intro}
Markov Decision Process (MDP) is an important framework modeling many applications in dynamic pricing, stochastic optimization and decision making \citep{Bertsekas, Puterman}. 
A stationary infinite horizon MDP is described by a set of states $\X$, sets of actions $\A_{s}$ for each state $s \in \X$, a transition kernel $\bm P$ which gives transition probabilities $\bm{P}_{sa} \in \R^{|\X|}_{+}$ for all state-action pair $(s,a)$, some rewards $r_{sa}$ for each state-action pair $(s,a)$ and a discount factor $ \lambda \in (0,1)$. 
A policy $\pi$ maps, for each period $t \in \N$, a state-action history up to time $t$ $(s_{0},a_{0},s_{1},a_{1},...,s_{t})$ to a probability distribution over the set of actions $\A_{s_{t}}$. 
Two important classes of policies are the set of {\em Markovian} policies and the set of {\em stationary} policies. A policy is called Markovian if it only depends on the current state $s_{t}$ and not on the complete history up to period $t$ $(s_{0},a_{0},s_{1},a_{1},...,s_{t})$.  A Markovian policy is called stationary if it does not depend on time.
We call $\Pi^{G}$ the set of all policies and $\Pi$ the set of stationary policies. The goal is to find a policy $\pi$ that maximizes the infinite horizon discounted expected reward 
 \begin{equation}\label{eq-expreward}
 R(\pi,\bm{P})=\E^{\pi, \bm{P}} \left[ \sum_{t=0}^{\infty} \lambda^{t}r_{s_{t}a_{t}} \; \bigg| \; s_{0} \sim \bm{p}_{0} \right]
 \end{equation}
{\color{black}where} $s_{t}$ is the state at period $t \in \N$ and $a_{t}$ is the action chosen at period $t$ following the probability distribution $(\pi_{s_{t}a})_{a} \in \R_{+}^{|\A|}$. The vector $\bm{p}_{0} \in \R_{+}^{|\X|}$ is a given initial probability distribution over the set of states $\X$. 
Without loss of generality, we assume $\A_{s}=\A$ for all states $s$. We also assume that the rewards are non-negative; otherwise, we can add the same large constant to all rewards,  without changing the optimal policy. We assume that the set of states and the set of actions are finite, with $|\X|=S$ and $|\A|=A$.
 
An optimal policy can be found in the set of stationary deterministic policies, which is potentially attractive from an implementation point of view. Efficient algorithms have been studied including policy iteration, value iteration and linear programming~\citep{Puterman,Bertsekas}. In particular, MDPs have been used extensively in many healthcare applications including early detection, prevention, screening and treatment of diseases. MDPs are particularly efficient to analyze chronic diseases and decisions that are made sequentially over time in a stochastic environment, see for instance the recent works on  cardiovascular controls for patients with Type 2 diabetes \citep{MDP-steimle} and determining the optimal stopping time for medical treatments \citep{cheng2019}.

However, in most applications, we only estimate model parameters from noisy observational data. Consequently, the decision maker may  recommend a policy that performs poorly with respect to the true parameters, since the optimal policy for the nominal parameters could potentially be highly sensitive to even small perturbations in the problem parameters and lead to highly suboptimal outcomes. 

We consider a robust approach to address the uncertainty in the transition probabilities where we model the uncertainty in $\bm{P}$ as an adversarial selection from a compact convex \textit{uncertainty set} $\PP$. 
Our goal is to find a policy that maximizes the worst-case expected reward over the choices of $\bm{P}$ in the uncertainty set $\PP$, i.e., our goal is to solve
\begin{equation}\label{eq:robust-mdp-problem}
z^{*}=\max_{\pi \in \Pi^{G}} \; \min_{\bm{P} \in \PP} \; R(\pi,\bm{P}).
\end{equation}
We refer to this as the \textit{policy improvement problem}, (PI), following the literature \citep{Kuhn}.  An important sub-problem of the above problem is to compute the worst-case reward of a given policy $\pi$,
\[z(\pi)=\min_{\bm{P} \in \PP} \; R(\pi,\bm{P}).\]
We refer to this as the \textit{policy evaluation problem} (PE).

The robust optimization approach has been specifically considered to address parameters uncertainty in MDPs, first in 1973  by~\cite{Satia},
and more recently in~\cite{Iyengar,Nilim,Kuhn}, and~\cite{Xu-Mannor}. In particular,
\cite{Iyengar} and~\cite{Nilim} consider a robust MDP where they model the uncertainty in transition probabilities using a {\em rectangular} uncertainty set, where transition probability $\bm{P}_{sa}$ for each state-action pair $(s,a)$ can be selected from a set $\PP_{sa} \subseteq \R_{+}^{|\X|}$, unrelated to transition probabilities out of other state-action pair, i.e., 
$ \PP = \underset{(s,a) \in \X \times \A}{\times} \;  \PP_{sa}, \text{ where } \; \PP_{sa} \subseteq \R^{|\X|}_{+}.$
They refer to this as an $(s,a)$-\textit{rectangular} uncertainty set
and show that for  such uncertainty sets, there is an optimal robust policy that is stationary  and deterministic and one can efficiently compute it using {\em robust value iteration}.

An important generalization of the $(s,a)$-rectangular uncertainty set is the \textit{row-wise} or $s$-\textit{rectangular} uncertainty set, introduced in~\cite{Epstein} and extensively studied in~\cite{Kuhn}. Here, the transition probabilities $\bm{P}_{s} = (P_{sas'})_{as'} \in \R_{+}^{|\A| \times |\X|}$ corresponding to different states are unrelated and the uncertainty set $\PP$ is given as:
 $\PP = \underset{s \in \X}{\times} \; \PP_{s}, \text{ where } \; \PP_{s} \subseteq \R^{|\X| \times |\A|}_{+}. $
\cite{Kuhn} show that for an $s$-rectangular uncertainty set, an optimal robust policy can be computed efficiently using a robust value iteration algorithm. An optimal robust policy can be chosen from the class of stationary policies but {\color{black}it} is not necessarily deterministic, which is in contrast with the case of classical MDPs or robust MDPs with $(s,a)$-rectangular uncertainty set, where there always exists an optimal policy that is deterministic. Importantly, \cite{Kuhn} also show that the problem of computing the optimal robust policy is strongly NP-Hard for general uncertainty set.

While the $(s,a)$-rectangularity and the $s$-rectangularity assumptions allow to design efficient algorithms to compute the optimal robust policy, such rectangular uncertainty sets are quite conservative in modeling uncertainty in {\color{black}many} applications. In particular, the rectangularity assumption allows the adversary to select the transition probabilities across different states unrelated to transitions in the other states. This is potentially very conservative especially if the transition probabilities depend on a common set of underlying factors.
{\color{black} For instance,  Markov model may be used to describe the health evolution of a patient (see \cite{mdp-med-1} for some examples). The state $s \in \X$ represents the health condition, and the action $a \in \A$ represents the treatment. The transitions $\bm{P}_{sa}$ represent the dynamics of the health evolution across different patient conditions, given a treatment.  In this healthcare setting,  
 common underlying factors, such as genetics, demographics and/or physiologic characteristics of certain diseases, could dictate the health evolution of patients across different health states.  Therefore,  the error in the estimation of one transition $\bm{P}_{sa}$ should be related to the errors in the estimations of the other transitions. As such, it is natural to consider a model of uncertainty where the transition probabilities of different health states are correlated.}

To address this issue, we consider a {\em factor uncertainty model}, where all transition probabilities depend on a small number of underlying factors. 
In particular, we consider $r$ factors, $\bm{w}_{1},..., \bm{w}_{r} \in \R^{\X}$ such that each transition probability $\bm{P}_{sa}$ is a linear (convex) combination of these factors. We refer to $\bm{W}  = (\bm{w}_{1},..., \bm{w}_{r})$ as the {\em factor matrix} and the model as {\em factor matrix uncertainty set}.  
The rank of the Markov chain induced by a stationary policy $\pi$ and a transition kernel $\bm{P}$ is (at most) $r$. \cite{Goh} also consider this model of uncertainty in the context of a healthcare application for modeling evolution of the health of the patient and give an algorithm to evaluate the worst-case expected reward of a policy under certain assumptions. However, they do not consider the problem of computing an optimal robust policy in this model.
 
\vspace{2mm}
\noindent \textbf{Our Contributions.}
Our main contributions are the following:

\vspace{2mm}
\noindent {\textbf {Generality of our model.}} We show that factor matrix uncertainty sets are very general and can model any uncertainty set. Therefore, we need an additional assumption to compute an optimal robust policy. In particular, we assume that each factor, $\bm{w}_i, i=1,\ldots,r$ can be chosen from an uncertainty set unrelated to other factors. Even with this assumption, we show that factor matrix uncertainty sets are general enough to model $(s,a)$-rectangularity in the special case where $r=S \cdot A$.

\vspace{2mm}
\noindent {\textbf{Min-max duality.}} 
We prove a structural min-max duality result:
\[ \max_{\pi \in \Pi^{G}} \; \min_{\bm{P} \in \PP} \; R(\pi,\bm{P}) = \min_{\bm{P} \in \PP}\;  \max_{\pi \in \Pi^{G}} \; R(\pi,\bm{P}).\]
In fact, we prove a stronger result: the same pair $(\pi^{*},\bm{P}^{*})$ is optimal on both sides of the strong duality equality. 
This implies that the optimal robust policy $\pi^{*}$ is an optimal policy for $\bm{P}^{*}$. Therefore, an optimal robust policy can be chosen stationary  and deterministic, which is useful for implementability and interpretability.

\vspace{2mm}
\noindent {\textbf{Optimal robust policy.}}
We give an efficient algorithm to compute an optimal robust policy. To do so, we first show that the evaluation of the worst-case of a policy can be reformulated as an alternate MDP. We then show that the problem of maximizing the worst-case reward can be reformulated as a coupled MDP, where the decision maker is playing a zero-sum game against an adversary. Computing an optimal robust policy then reduces to finding the fixed-point of a contraction. This yields an efficient algorithm for finding an optimal robust policy, using robust value iteration.
To the best of our knowledge, this is the first example of an uncertainty set where the transition probabilities across different states are related and still one can compute the optimal robust policy. 

\vspace{2mm}
\noindent {\textbf{Maximum principle and Blackwell optimality}}. We show that certain important structural properties of classical MDP also hold for the optimal robust policy for factor matrix uncertainty sets. In particular, we present the \textit{robust maximum principle}, which states that the worst-case value function of an optimal robust policy is component-wise higher than the worst-case value function of any other policy. Moreover, we prove the robust counterpart of \textit{Blackwell optimality}, which states that there exists a pair $(\pi^{*},\bm{P}^{*})$ that remains optimal for the policy improvement problem for all discount factor sufficiently close to $1,$ and $\pi^{*}$ can be chosen deterministic.

\vspace{2mm}
\noindent {\textbf{Numerical experiments.}}
We detail the computation of the factor matrix from the estimated nominal kernel and we highlight the role of nonnegative matrix factorization. We show on two MDP examples that the performances of the optimal nominal policy can significantly deteriorate for small variations of the parameters and we compare the performances of robust policies related to factor matrix and to $s$-rectangular uncertainty sets. We show that our optimal robust policy both improves the worst-case reward and has better nominal reward than the robust policy related to the $s$-rectangular set. Our robust policy also has better empirical performances than the robust policy of the $s$-rectangular uncertainty set. Our results suggest that the factor matrix uncertainty set is a less conservative model to handle uncertainty in parameters than the $s$-rectangular uncertainty set.

\section{Factor matrix uncertainty set.}
\label{sec:model}
 We consider an uncertainty set $\PP \subseteq \R^{S \times A \times S}_{+}$ of the form
\begin{equation}\label{eq:fmus}
\PP = \left\lbrace \left( \sum_{i=1}^{r} u^{i}_{sa}w_{i,s'} \right)_{sas'}  \bigg| \;    \bm{W} = (\bm{w}_{1},..., \bm{w}_{r}) \in \W \subseteq \R^{S \times r} \right\rbrace
\end{equation} 
where $ \; \bm{u}_{1},..., \bm{u}_{S}$ are fixed and known in $\R_{+}^{r \times A}$ and $\W$ is a convex, compact subset of $\R_{+}^{S \times r}$ such that: 
\begin{align}\label{eq:uisa-distrib}
\sum_{i=1}^{r} u^{i}_{sa}  = 1, \; \forall \; (s,a) \in \X \times \A, 
\sum_{s'=1}^{S} w_{i,s'}  = 1, \; \forall \; i \in [r],
\end{align}
 where we write $[r] = \{1,...,r\}$. Note that here $r$ is an integer and is not related to the rewards $r_{sa}$ for $(s,a) \in \X \times \A$.
We refer to the above uncertainty set as \textit{factor matrix uncertainty set}. Such a model allows to capture correlations across transitions probabilities, unlike $s$- and $(s,a)$-rectangular uncertainty sets. Each transition vector $\bm{P}_{sa} \in \R^{S}_{+}$ is a convex combination of some common underlying \textit{factors} $\bm{w}_{1},...,\bm{w}_{r}$ with fixed
\textit{coefficients} $u^{1}_{sa}, ..., u^{r}_{sa}$, i.e., $
\bm{P}_{sa} = \sum_{i=1}^{r} u^{i}_{sa}\bm{w}_{i},$
where each of the factor $\bm{w}_{i}$ is a probability distribution over the next state in $\X$. Since for all $(s,a) \in \X \times \A$ the vectors $\bm{P}_{sa}$ are convex combination of the same factors $\bm{w}_{1},...,\bm{w}_{r}$, 
 this class of uncertainty sets can model coupled transitions. We detail the computation of the factors $\bm{w}_{1}, ..., \bm{w}_{r}$ and the coefficients $\bm{u}_{1}, ..., \bm{u}_{r}$ in Section \ref{sec:simu} where we highlight the role of nonnegative matrix factorization.
Note that this model has been considered in \cite{Goh} in the context of robust MDP.

\vspace{2mm}
\noindent
{\color{black}
We start by proving that factor matrix uncertainty set are very general. In particular,  we have the following proposition.
\begin{proposition}\label{prop:generality-r-rec}
Let $\PP \subset \R^{S \times A \times S}$ be any uncertainty set (not necessarily $s$-rectangular or $(s,a)$-rectangular). Then $\PP$ can be reformulated as a factor matrix uncertainty set.
\end{proposition}
\proof{Proof.}
Given an uncertainty set $\PP$,  we choose $r=S \cdot A, \W = \PP$ and $u_{sa}^{(\hat{s},\hat{a})} = 1$ if $(s,a)=(\hat{s},\hat{a})$, and $0$ otherwise.
In this case, it is clear that for all $(s,a) \in \X \times \A$,
\[ \sum_{i=1}^{r} u^{i}_{sa} = \sum_{(\hat{s},\hat{a}) \in \X \times \A} u_{sa}^{(\hat{s},\hat{a})} = 1,\]
and that for all $ i \in [r]$, i.e., for all $(\hat{s},\hat{a}) \in \X \times \A$,
\[ \sum_{s'=1}^{S} w_{i,s'} = \sum_{s'=1}^{S} P_{\hat{s}\hat{a}s'} = 1.\]
Finally,   for all $(s,a) \in \X \times \A$,
\[  \sum_{i=1}^{r} u^{i}_{sa}\bm{w}_{i} =  \sum_{(\hat{s},\hat{a}) \in \X \times \A} u_{sa}^{(\hat{s},\hat{a})} \bm{w}_{(\hat{s},\hat{a})} = \bm{w}_{(s,a)} = \bm{P}_{sa}.\]
\hfill \qed
\endproof
Note that in the proof of the previous proposition, we choose $r= S \cdot A$ and then enumerate the vectors $\left( \bm{P}_{sa} \right)_{(s,a) \in \X \times \A}$ with the factors $\bm{w}_{1}, ...., \bm{w}_{r}$. }

Since factor matrix uncertainty sets are able to model any uncertainty set, the problem of finding the worst-case transition $\bm{P}$ for a given policy $\pi \in \Pi^{G}$ is intractable. Indeed,\cite{Kuhn} and \cite{Goh} show that for any fixed policy $\pi \in \Pi^{G}$, general uncertainty set $\PP$ and scalar $\gamma$, it is NP-complete to decide whether $\min_{\bm{P} \in \PP} \; R(\pi,\bm{P}) \geq \gamma.$

The intractability of the policy evaluation problem indicates that the policy improvement problem is also intractable.
In view of the intractability for general factor matrix uncertainty set, we make the following additional assumption on the set $\W$.
\begin{assumption}[$r$-rectangularity]\label{ass:r-rec}
$\W = \W^{1} \times ... \times \W^{r}, \text{ where } \W^{1}, ..., \W^{r} \subset \R^{S}_{+}. $
\end{assumption}
We refer to this property as \textit{$r$-rectangularity}. Factors $\bm{w}_{1},...,\bm{w}_{r}$ are said to be \textit{unrelated}, because each vector $\bm{w}_{i}$ can be selected in each set $\W^{i}$ unrelated to $\bm{w}_{j}$, for $j \neq i$. 
For any state-action pair $(s,a)$, the factors $\bm{w}_{1},...,\bm{w}_{r}$ are combined to form the transition kernel $\bm{P}$. Therefore, $r$-rectangular uncertainty sets also model relations between the transition probabilities across different states. 
{\color{black}
We detail the relations between $r$-, $s$- and $(s,a)$-rectangularity in the following proposition.
\begin{proposition}\label{prop:relationship-s-r-sa-rec}
\begin{enumerate}
\item An $(s,a)$-rectangular uncertainty set is also $r$-rectangular.
\item There exists an $r$-rectangular uncertainty set that is not $s$-rectangular.
\item There exists an $r$-rectangular uncertainty set that is not $(s,a)$-rectangular.
\item There exists an $s$-rectangular uncertainty set that is not $r$-rectangular.
\end{enumerate}
\end{proposition}
\proof{Proof.}
\begin{enumerate}
\item From the proof of Proposition \ref{prop:generality-r-rec}, any uncertainty set $\PP$ can be formulated as a factor matrix uncertainty set with $r=S \cdot A, \W = \PP$, and
\[\bm{P}_{sa} =  \bm{w}_{(s,a)} = \sum_{(\hat{s},\hat{a}) \in \X \times \A} u_{sa}^{(\hat{s},\hat{a})} \bm{w}_{(\hat{s},\hat{a})}, \forall \; (s,a) \in \X \times \A,\]
for $u_{sa}^{(\hat{s},\hat{a})} = 1$ if $(s,a)=(\hat{s},\hat{a})$, and $0$ otherwise.

If $\PP$ is an $(s,a)$-rectangular uncertainty set, then the vectors $\left(\bm{P}_{sa}\right)_{(s,a) \in \X \times \A}$ are unrelated and so do the factors $\bm{w}_{1},...,\bm{w}_{r}=\bm{w}_{(1,1)},...,\bm{w}_{(S,A)}$. Therefore $\PP$ is an $r$-rectangular set.
\item This is because $r$-rectangularity can model correlations in transitions across different states.  For example, consider a robust MDP where there are only two states $s_{1}$ and $s_{2}$, one action $a$ and $r=1$. In such a case, there exists a set $\W \subseteq \R_{+}^{S}$ such that
\[\bm{P}_{s_{1}a}=\bm{P}_{s_{2}a}=\bm{w} \in \W,\] and therefore the uncertainty set $\PP$ is $r$-rectangular. In particular,
\begin{align*}
\PP &  = \{(\bm{P}_{s_{1}a},\bm{P}_{s_{2}a}) \; | \; \bm{P}_{s_{1}a}=\bm{P}_{s_{2}a}, \bm{P}_{s_{1}a} \in \W \} \\
&  = \{(\bm{w},\bm{w}) \; | \; \bm{w} \in \W \}.
\end{align*}
However, the set $\PP$ is not $s$-rectangular, because the smallest $s$-rectangular set containing $\PP$ is
\[\{(\bm{w_{1}},\bm{w_{2}}) \; | \; \bm{w_{1}},\bm{w_{2}} \in \W \} = \W \times \W,\]
and the set $\W \times \W$ is different from $\PP.$
\item \tb{The proof from the previous point also applies for $(s,a)$-rectangular uncertainty set, since in the previous example we choose $a=1$.}
\item We defer the proof of this statement to Section \ref{sec:duality}, Remark \ref{rem:s-rec-not-r-rec}.
\end{enumerate}
\hfill \qed
\endproof
From Proposition \ref{prop:relationship-s-r-sa-rec}, we see that $r$-rectangularity is a generalization of $(s,a)$-rectangularity, different from $s$-rectangularity, which also generalizes $(s,a)$-rectangularity. The main difference is that $r$-rectangularity can model coupled transitions across states, but requires a common \textit{linear} (convex) dependency on some underlying vectors $\bm{w}_{1}, ..., \bm{w}_{r}$. For $s$-rectangularity, the transitions are uncoupled across states, but the vectors $(\bm{P}_{sa})_{a \in \A}$ (for any $s \in \X$) do not need to linearly rely on some underlying factors.

\tb{
We also discuss the relation between $r$-rectangularity and $k$-rectangularity \citep{Xu-Mannor}.  
In particular, $k$-rectangularity specifically focuses on MDPs with \textit{finite horizon} $h \in \N$,  with stationary transition kernel $\bm{P} \in \R^{S \times A \times S}$ for all periods $0,...,h$.  In short, $k$-rectangularity assumes that for any subset of states $\X' \subset \X$, and conditioned on parameter realization at all $s \in \X'$, the set of possible values for the other transitions $\bm{P}_{s} \in \R^{A \times S}$ (for $ s \in \X \setminus \X'$) belongs to at most $k$ different possible sets. Therefore, $k$-rectangularity can also account for non-rectangular parameter deviations. However, it requires to augment the state space from $\X$ to $\X \times [k]$, i.e., at period $t$  the decision maker must be able to observe the values of the kernels $\bm{P}_{s_{t'}}$ at every state $s_{t'}$ at every period $t' \leq t$, and keep counts of the number of different kernels observed so far. This is in stark contrast with $r$-rectangular uncertainty set, where we consider infinite horizon discounted rewards, and the decision maker does not need to observe the entire transition kernel $\bm{P}_{s_{t}}$ at each period $t$.}
   
\tb{In the rest of the paper, we will use the following assumptions for the structural properties of the uncertainty sets $\W^{1}, ..., \W^{r}$.}
\begin{assumption}\label{ass:tractable-u-set}
\tb{The sets $\W^{1}, ..., \W^{r}$ are convex compact. Moreover, for any $i \in [r]$, for any cost vector $\bm{c} \in \R^{S}$, we can find an $\epsilon$-optimal solution of
\[ \min_{\bm{w}_{i} \in \W^{i}} \bm{c}^{\top}\bm{w}_{i}\]
 in $O \left( comp\left(\W^{i}\right) \log\left(\epsilon^{-1}\right)\right)$ arithmetic operations, where $comp(\W^{i}) \in \R$ depends on the structure of the uncertainty set $\W^{i}$.}
\end{assumption}
%
\tb{Assumption \ref{ass:tractable-u-set} is very general. In particular, it is satisfied for many classical uncertainty sets, such as polyhedral uncertainty sets (e.g.,  based on $\ell_{1}$- and weighted $\ell_{1}$-norms~\citep{Iyengar,Ho}, $\ell_{\infty}$-norm~\citep{givan1997bounded}, and budget of uncertainty sets~\citep{BS2004}) or ellipsoidal-based uncertainty sets~\citep{ben2000robust,Kuhn}. We refer to Section 4.6 in \cite{BenTal-Nemirovski} for the precise value of $comp(\W^{i})$ in these cases.}

We now consider the properties of optimal policies. For $r$-rectangular uncertainty sets, there exists an optimal robust policy that is stationary.  In particular, we have the following proposition.
\begin{proposition}\label{prop:stat-mrkv} 
Let $\PP$ be an $r$-rectangular uncertainty set. \tb{Under Assumption \ref{ass:tractable-u-set},} there exists a stationary policy that is a solution to the policy improvement problem.
\end{proposition}
We present the proof in Appendix \ref{app:stat-mrkv}. In view of Proposition \ref{prop:stat-mrkv}, in the rest of the paper we focus on policies in the set $\Pi$ of stationary polices (possibly randomized).
\section{Policy evaluation for $r$-rectangular uncertainty set.}\label{sec:PE}
In this section we consider the policy evaluation problem, where the goal is to compute the worst-case transition kernel of a policy. \cite{Goh} give an algorithm for the policy evaluation problem for $r$-rectangular factor matrix uncertainty sets. To compute an optimal robust policy in Section \ref{sec:PI}, we present in this section an alternate algorithm for the policy evaluation problem, which provides structural insights on the solutions of the policy improvement problem. 
\subsection{Algorithm for policy evaluation.}
We show that the policy evaluation problem can be written as an alternate MDP with $r$ states and $S$ actions, and $\W$ as the set of policies. This alternate MDP is played by the \textit{adversary}.
 Let us introduce some notations to formalize our intuition. We fix a policy $\pi \in \Pi$ and we write
$\bm{T}_{\pi} = \left(\sum_{a \in \A} \pi_{sa}u^{i}_{sa}\right)_{(s,i) \in \X \times [r]} \in \R_{+}^{S \times r}.$

\vspace{2mm}
\noindent {\textbf {Adversarial MDP}}.
The policy evaluation problem can be reformulated as an MDP with state set $[r]$, action set $\X$ and policy set $\W$. This MDP is played by the adversary. The adversary starts at period $t=0$ with an initial reward $\bm{p}_{0}^{\top}\bm{r_{\pi}}$ and the initial distribution $\bm{p}_{0}^{\top}\bm{T}_{\pi}$ over the set of states $[r]$. When the current state is $i \in [r]$, the adversary picks action $s$ with probability $w_{i,s}$. For $i,j \in [r]$ and $ s \in \X$, the transition probability and the reward are given by
\begin{equation}\label{eq:adversarial-MDP}
Prob \left( i \underset{action \; s}{\rightarrow} \; j \right) =\sum_{a \in \A} \pi_{sa}u^{j}_{sa}, \; Reward \left(i, \; action \; s\right) = \sum_{a \in \A} \pi_{sa}r_{sa}.
\end{equation}
 It is worth noting that the transition probability only depends on the chosen action $s$ and the arriving state $j$ but not on the current state $i$. The reward only depends on the chosen action $s$ and not on the current state $i$.   Additionally,  the policy of the adversary may be randomized, as the factors $\bm{w}_{i}$ may not be deterministic. \tb{Also, note that the authors in \cite{ho2020partial} construct a \textit{robust evaluation MDP}, with the goal to reformulate the policy evaluation problem for $s$-rectangular uncertainty set as an MDP. The main differences with our adversarial MDP are  the sets of states and actions, and the constraints on the set of feasible policies.} We have the following proposition.
 \begin{proposition} Let $\bm{\beta} \in \R^{r}$ be the value function of the adversary choosing factors $\bm{w}_{1} \in \W^{1}, ..., \bm{w}_{r} \in \W^{r}$. Then $\bm{\beta}$ satisfies the Bellman recursion:
 \begin{equation}\label{eq:bellman-recursion-adversary}
 \beta_{i} = \bm{w}_{i}^{\top}(\bm{r}_{\pi} + \lambda \cdot  \bm{T}_{\pi} \bm{\beta}), \forall \; i \in [r].
 \end{equation}
 \end{proposition} 
 \proof{Proof.}
The Bellman recursion for an MDP $(\X,\A,\lambda,\bm{P},\bm{r})$ is given by
\[ v_{s}^{\pi} = \sum_{a \in \A} \pi_{sa} \left( r_{sa} + \lambda \bm{P}_{sa}^{\top} \bm{v}^{\pi} \right), \forall \; s \in \X, \]
where $\bm{v}^{\pi}$ is the value function for the policy $\pi$ (\cite{Puterman}, Section 6.2.1).  We obtain \eqref{eq:bellman-recursion-adversary} by writing the Bellman recursion for the Adversarial MDP, where the set of states is $[r]$, the set of actions is $\X$, the set of  policies is $\W$,  and the transition rates and rewards are given by \eqref{eq:adversarial-MDP}.
\hfill \qed
 \endproof
 
The $r$-rectangularity assumption enables us to develop an iterative algorithm for the policy evaluation problem.
In particular, following the interpretation of the policy evaluation problem as an alternate MDP,  we have the following theorem.
We present a detailed proof in Appendix \ref{app:proof-PE}.
\begin{theorem}\label{th:PE-tractable} Let $\PP$ be an $r$-rectangular uncertainty set and $\pi$ be a stationary policy.
{\color{black}
Consider the following value iteration algorithm for the adversarial MDP, given by Algorithm \ref{alg:VI-PE}:
\begin{equation}\label{alg:VI-PE}\tag{VI-PE}
\bm{\beta}^{0} = \bm{0}, \beta^{k+1}_{i} =  \min_{\bm{w}_{i} \in \W^{i}} \; \bm{w}_{i}^{\top}\left( \bm{r}_{\pi} + \lambda \cdot \bm{T}_{\pi}\bm{\beta}^{k} \right), \forall \; i \in [r], \forall \; k \geq 0.
\end{equation}
}
\tb{Under Assumption \ref{ass:tractable-u-set}, Algorithm \ref{alg:VI-PE} returns an $\epsilon$-optimal solution to the policy evaluation problem in $O\left(rSA + rS \log(\epsilon^{-1}) + \sum_{i=1}^{r} comp(\W^{i}) \log^{2}\left(\epsilon^{-1}\right)\right)$ arithmetic operation.}
\end{theorem}
 
\vspace{2mm}
\noindent \textbf{Role of $r$-rectangularity.}
In the classical MDP framework, one assumes that the decision maker can independently choose the distributions $\pi_{s}$ across different states. This is because the set of stationary policies $\Pi$ is itself a Cartesian product:
\[\Pi = \left\lbrace (\pi_{sa})_{s,a} \in \R_{+}^{S \times A} \; \bigg| \; \forall s \in \X, \sum_{a \in \A} \pi_{sa} =1 \right\rbrace = \underset{s \in \X}{\times} \; \left\lbrace \bm{p} \in \R_{+}^{A} \; \bigg| \; \sum_{a \in \A} p_{a} =1 \right\rbrace. \]
Using the rectangularity of the policy set, one derives a fixed-point equation for the value function of the MDP from the classical Bellman Equation:  $v^{*}_{s} = \max_{a \in \A} \left\{ r_{sa} + \lambda \cdot \bm{P}_{sa}^{\top}\bm{v}^{*} \right\}, \forall \; s \in \X$ where $\bm{v}^{*}$ is the value function of the optimal policy {\color{black} (\cite{Puterman}, Section 6.2.1).}
This is the basis of the analysis of value iteration, policy iteration and linear programming algorithms for MDPs.

If the set $\PP$ is not $r$-rectangular, i.e., if there are some constraints across the factors $\bm{w}_{1},...,\bm{w}_{r}$, the adversary can not optimize independently over each component of the vector $\bm{\beta}$ since the same factor $\bm{w}_{i}$ can be involved in different components of the Bellman Equation \eqref{eq:bellman-recursion-adversary}. In particular, the factors $\bm{w}^{*}_{1},...,\bm{w}^{*}_{r}$ that attain the minima in \eqref{eq:bellman-recursion-adversary} might not be feasible in $\W$. We provide a counter-example for this in Appendix \ref{app:proof-PE}. However, when the uncertainty set $\PP$ is $r$-rectangular, the set $\W$ is a Cartesian product and the adversary can optimize independently over each of the sets $\W^{1},...,\W^{r}$ and recover a feasible solution in $\W = \W^{1} \times ... \times \W^{r}$, as detailed in Algorithm \ref{alg:VI-PE}.

\subsection{A reformulation of the policy evaluation problem}\label{sec:reform-PE-general}
We introduce here a reformulation of $z(\pi)$, which is useful to analyze the structure of the set of optimal robust policies in the next section.
 The proof of the next lemma relies on a fixed-point theorem and is detailed in Appendix \ref{app:PE:reform-general}.
\begin{lemma}\label{lem:PEfinal} \tb{Under Assumption \ref{ass:tractable-u-set}, the policy evaluation problem can be reformulated as follows.
\begin{align}\label{eq:PE-lemPEfinal}
z(\pi) = \; \max \;  & \bm{p}_{0}^{\top}(\bm{r}_{\pi} + \lambda \cdot \bm{T}_{\pi} \bm{\beta}) \\
& \beta_{i} \leq \min_{\bm{w}_{i} \in \W^{i}} \; \bm{w}_{i}^{\top}\left( \bm{r}_{\pi} + \lambda \cdot \bm{T}_{\pi}\bm{\beta} \right), \; \forall \; i \in [r], \\
& \bm{\beta} \in \R^{r}.
\end{align}
Moreover, an optimal solution $\bm{\beta}^{\pi}$ can be chosen such that for any $i \in [r]$, 
\[\beta_{i}^{\pi} = \min_{\bm{w}_{i} \in \W^{i}} \; \bm{w}_{i}^{\top}\left( \bm{r}_{\pi} + \lambda \cdot \bm{T}_{\pi}\bm{\beta}^{\pi} \right) .\]}
\end{lemma}

\section{Min-max duality.}\label{sec:duality}
In this section, we analyze the structure of the set of optimal robust policies. In particular, we present our min-max duality result.
Using reformulation \eqref{eq:PE-lemPEfinal}, the policy improvement problem $z^{*}$ becomes
\begin{align}\label{eq:PIfinal}
z^{*} = \max_{\pi \in \Pi} \;  z(\pi) = \; \max \;  & \bm{p}_{0}^{\top}(\bm{r}_{\pi} + \lambda \cdot \bm{T}_{\pi} \bm{\beta}) \\
& \beta_{i} \leq \min_{\bm{w}_{i} \in \W^{i}} \; \bm{w}_{i}^{\top}\left( \bm{r}_{\pi} + \lambda \cdot \bm{T}_{\pi}\bm{\beta} \right), \; \forall \; i \in [r], \\
& \bm{\beta} \in \R^{r}, \pi \in \Pi.
\end{align} 
In Proposition \ref{prop:stat-mrkv}, we have shown that there is a robust optimal policy that is stationary. In the following lemma, we show that a robust optimal policy can be chosen stationary  and deterministic; this is potentially attractive from an application standpoint.
\begin{lemma}\label{lem:determ-poli}
Let $\PP$ be an $r$-rectangular uncertainty set.
\tb{Under Assumption \ref{ass:tractable-u-set},} there exists a stationary  and deterministic policy solution of the policy improvement problem.
\end{lemma}
\proof{Proof.}
\tb{
Consider $\pi^{*}$ an optimal robust policy and $(\pi^*, \bm{\beta}^*)$ an optimal solution of  
\eqref{eq:PIfinal}.
Consider the following policy $\tilde{\pi}$ where for all $s \in \X$,\[\tilde{\pi}_{s} \in \arg  \max_{\pi_{s} \in \Delta} \; \sum_{a \in \A} \pi_{sa} \left(  r_{sa} + \lambda \cdot  \sum_{i=1}^{r}  u^{i}_{sa}\bm{\beta}^{*}_{i} \right),\]
where we write $\Delta$ for the simplex over $\A$: $\Delta = \{ \bm{p} \in \R^{\A}_{+} \; | \; \sum_{a \in \A} p_{a}=1\}$.}
\tb{
The policy $\tilde{\pi}$ can be chosen deterministic, because for each $s \in \X$ the distribution $\tilde{\pi}_{s}$ is a solution of a linear program over the simplex $\Delta$, and the extreme points of $\Delta$ are the distributions over $\A$ with exactly one non-zero coefficient, and this coefficient is equal to $1$.
By definition, the deterministic policy $\tilde{\pi}$ has an objective value in \eqref{eq:PIfinal} at least as high as the objective value of $\pi^{*}$.
 Moreover, $(\tilde{\pi},\bm{\beta}^{*})$ is still feasible in \eqref{eq:PIfinal}.}
 \tb{Indeed, because $\left(\pi^{*},\bm{\beta}^{*}\right)$ is feasible, we have
\begin{equation}\label{eq:feasibility-pi-start}
\bm{\beta}^{*} \leq \bm{r}_{\pi^{*}} + \lambda \cdot \bm{T}_{\pi^{*}}\bm{\beta}^{*}.
\end{equation}
By definition of $\tilde{\pi}$, we have
\begin{equation}\label{eq:comp-pi-star-pi-tilde}
\bm{r}_{\pi^{*}} + \lambda \cdot \bm{T}_{\pi^{*}}\bm{\beta}^{*} \leq  \bm{r}_{\tilde{\pi}} + \lambda \cdot \bm{T}_{\tilde{\pi}}\bm{\beta}^{*}.
\end{equation}
Therefore, by combining \eqref{eq:feasibility-pi-start} and \eqref{eq:comp-pi-star-pi-tilde}, we obtain 
\[ \bm{\beta}^{*} \leq \bm{r}_{\tilde{\pi}} + \lambda \cdot \bm{T}_{\tilde{\pi}}\bm{\beta}^{*}.\]
This proves that $\tilde{\pi}$ is feasible.
Therefore, there exists a stationary and deterministic policy solution to the policy improvement problem.}
\hfill \qed

 \tb{Since each transition kernel $\bm{P} \in \PP$ is fully determined by a factor matrix $\bm{W} \in \W$, for the rest of the paper we write $R(\pi,\bm{W})$, instead of $R(\pi,\bm{P})$, for the infinite horizon discounted expected reward, defined in \eqref{eq-expreward}. }
We will now prove our min-max duality result.
\begin{theorem}\label{th:duality} \tb{Under Assumption \ref{ass:tractable-u-set},} let $(\pi^{*}, \bm{W}^{*})$ be a solution to {\color{black}the robust MDP problem \eqref{eq:robust-mdp-problem} with $r$-rectangular uncertainty set,} with $\pi^{*}$ deterministic. Then
\begin{equation}\label{eq:argminargmax}
\bm{W}^{*} \in \arg \min_{\bm{W} \in \W} R(\pi^{*},\bm{W}) \text{ and } \pi^{*} \in \arg \max_{\pi \in \Pi} R(\pi,\bm{W}^{*}).
\end{equation}
Moreover, the following strong min-max duality holds.
\begin{equation}\label{eq:strong-min-max-duality}
\max_{\pi \in \Pi} \; \min_{\bm{W} \in \W} \; R(\pi,\bm{W}) = \min_{\bm{W} \in \W}\;  \max_{\pi \in \Pi} \; R(\pi,\bm{W}).
\end{equation}
\end{theorem}
We have shown in Section \ref{sec:PE} that the policy evaluation problem can be reformulated as an alternate MDP, played by the adversary. We also introduced $\bm{\beta} \in \R^{r}$ the value function for the adversary, defined by the Bellman Equation \eqref{eq:bellman-recursion-adversary}. To prove Theorem \ref{th:duality}, we need the following lemma that relates the value function $\bm{v}$ of the decision maker and the value function $\bm{\beta}$ of the adversary. 

\begin{lemma}\label{lem:Wv-beta}
Let $\pi \in \Pi$ and $\bm{W} \in \W$. Let $\bm{v}$ be the value function of the decision maker and $\bm{\beta}$ be the value function of the adversary. Then $
\bm{W}^{\top}\bm{v} = \bm{\beta}.$
\end{lemma}
\proof{Proof.}

\tb{By definition, the vector $\bm{\beta}$ is the unique solution of the Bellman recursion in the adversarial MDP:
\begin{equation}\label{eq:Bell-beta-3}
\beta_{i} = \bm{w}_{i}^{\top} \left(\bm{r}_{\pi} + \lambda \bm{T}_{\pi} \bm{\beta} \right), \forall \; i \in [r].
\end{equation}
Therefore, to prove that $\bm{W}^{\top}\bm{v}=\bm{\beta}$, we simply need to prove that $\bm{W}^{\top}\bm{v}$ satisfies 
\[ \left(\bm{W}^{\top}\bm{v} \right)_{i} = \bm{w}_{i}^{\top}\left(\bm{r}_{\pi} + \lambda \bm{T}_{\pi} \bm{W}^{\top}\bm{v} \right), \forall \; i \in [r].\]
Note that $\bm{W}^{\top}\bm{v} = (\bm{w}_{i}^{\top}\bm{v})_{i \in [r]}$.
We start from the Bellman Equation for the decision maker: the value function $\bm{v}$ is uniquely determined by 
\begin{equation}\label{eq:Bell-2}
v_{s} = \sum_{a \in \A} \pi_{sa}(r_{sa} + \lambda \cdot \sum_{i=1}^{r} u^{i}_{sa}\bm{w}_{i}^{\top}\bm{v}), \forall \; s \in \X.
\end{equation}
This equation can be written more concisely as
\begin{equation}\label{eq:Bell-new}
\bm{v} = \bm{r}_{\pi} + \lambda \bm{T}_{\pi} \bm{W}^{\top}\bm{v}.
\end{equation}
We multiply the equation in \eqref{eq:Bell-new} by $\bm{w}_{i}^{\top}$, for $i \in [r]$, and we obtain 
\[\bm{w}_{i}^{\top}\bm{v}=\bm{w}_{i}^{\top}(\bm{r}_{\pi} + \lambda \cdot \bm{T}_{\pi}\bm{W}^{\top}\bm{v}), \forall \; i \in [r].\]
Therefore, $\bm{W}^{\top}\bm{v}$ satisfies the Bellman Equation for the adversarial MDP \eqref{eq:Bell-beta-3}. Because $\bm{\beta}$ is the unique solution to \eqref{eq:Bell-beta-3}, we conclude that $\bm{W}^{\top}\bm{v} = \bm{\beta}.$
}
\hfill \qed

We are now ready to prove Theorem \ref{th:duality}.

\proof{Proof of Theorem \ref{th:duality}.}
When the decision maker chooses policy $\pi^{*}$ and the adversary chooses factor matrix $\bm{W}^{*}$, let $\bm{v}^{*}$ be the value function of the decision maker and $\bm{\beta}^{*}$ be the value function of the adversary. Let $\bm{P}^{*}$ be the transition kernel associated with the factor matrix $\bm{W}^{*}$. 
Since $\pi^{*}$ is deterministic, we write $a^{*}(s) \in \A$ the action chosen in each state $s$, uniquely determined by $ \forall \; s \in \X, \; \pi^{*}_{s a^{*}(s)} = 1.$

Our goal is to show $ \pi^{*} \in \arg \max_{\pi \in \Pi} R(\pi,\bm{W}^{*}).$
From the Bellman Equation, this is equivalent to showing that
$
 v^{*}_{s} = \max_{a \in \A} \left\{ r_{sa} + \lambda \cdot \bm{P}_{sa}^{* \; \top}\bm{v}^{*} \right\}, \forall \; s \in \X.
$
In the case of an $r$-rectangular uncertainty set, this is equivalent to proving that
\begin{equation}\label{eq:pistar:argmax}
v^{*}_{s} = \max_{a \in \A} \left\{ r_{sa} + \lambda \cdot \sum_{i=1}^{r} u^{i}_{sa}\bm{w}_{i}^{* \top}\bm{v}^{*} \right\}, \forall \; s \in \X.
\end{equation}
Since $\bm{v}^{*}$ is the value function of the policy $\pi^{*}$ when the adversary picks $\bm{W}^{*}$, it satisfies the Bellman Equation \eqref{eq:Bell-2}. For all $s \in \X$,
\begin{align}
v^{*}_{s} & \label{eq:opt1} =  r_{sa^{*}(s)} + \lambda \cdot \sum_{i=1}^{r} u^{i}_{sa^{*}(s)}\bm{w}_{i}^{*\top}\bm{v}^{*}
 \\
& \label{eq:opt2}   =  r_{sa^{*}(s)} + \lambda \cdot \sum_{i=1}^{r} u^{i}_{sa^{*}(s)}\beta_{i}^{*}
 \\
& \label{eq:opt4} = \max_{a \in \A} \left\{r_{sa} + \lambda \cdot \sum_{i=1}^{r} u^{i}_{sa}\beta_{i}^{*} \right\} \\
& \label{eq:opt5} =\max_{a \in \A} \left\{ r_{sa} + \lambda \cdot \sum_{i=1}^{r} u^{i}_{sa}\bm{w}_{i}^{*\top}\bm{v}^{*} \right\},
\end{align}
where \eqref{eq:opt1} follows from Bellman Equation on the deterministic policy $\pi^{*}$ and \eqref{eq:opt2} follows from Lemma \ref{lem:Wv-beta}.
 The key Equality \eqref{eq:opt4} follows from the choice of $\pi^{*}$ deterministic as in Lemma \ref{lem:determ-poli}. Finally, \eqref{eq:opt5} is a consequence of Lemma \ref{lem:Wv-beta}.
 
Following \eqref{eq:pistar:argmax}, we conclude that $\pi^{*} \in \arg \max_{\pi \in \Pi} \; R(\pi,\bm{W}^{*}).$
Note that classical weak duality holds:
\begin{equation}
\label{ineq:maxmin-minmax}
\max_{\pi \in \Pi} \; \min_{\bm{W} \in \W} \; R(\pi,\bm{W}) \leq \min_{\bm{W} \in \W}\;  \max_{\pi \in \Pi} \; R(\pi,\bm{W}).
\end{equation} 
But the pair $(\pi^{*},\bm{W}^{*})$ bridges this gap, because 
$
\bm{W}^{*} \in \arg\min_{\bm{W} \in \W} R(\pi^{*},\bm{W}) \text{ and } \pi^{*} \in \arg \max_{\pi \in \Pi} R(\pi,\bm{W}^{*}).
$
Therefore the two sides of \eqref{ineq:maxmin-minmax} are attained at the same pair $(\pi^{*},\bm{W}^{*})$ and we obtain our strong duality result.
\hfill \qed

\begin{remark}\label{rem:s-rec-not-r-rec}
\tb{
In Theorem~\ref{th:duality}, the fact that the optimal robust policy $\pi^{*}$ is the optimal nominal policy of the MDP where the adversary plays $\bm{W}^{*}$ can be seen as an {\em equilibrium} for the zero-sum game between the decision maker and the adversary. This result is in contrast with the case of $s$-rectangular uncertainty sets. Indeed, we provide in Appendix \ref{app:counter-example-equilibrium} an example of an MDP instance where the uncertainty set is $s$-rectangular and the optimal solutions of the left-hand side and the right-hand side of \eqref{eq:strong-min-max-duality} are different. This also proves the existence of an $s$-rectangular uncertainty set that can not be reformulated as an $r$-rectangular uncertainty set (see Proposition \ref{prop:relationship-s-r-sa-rec}).}
\end{remark}

\section{Policy improvement for $r$-rectangular uncertainty set.}\label{sec:PI}
We consider the policy improvement problem where we want to find a policy with the highest worst-case reward. We give an efficient algorithm to compute an optimal robust policy, assuming that the set $\PP$ is $r$-rectangular. 

We have shown in Lemma \ref{lem:determ-poli} that there exists a deterministic optimal robust policy. This motivates us to consider the following iterative algorithm, Algorithm \ref{alg:VI-PI}, that computes a deterministic policy in each iteration. In particular, in each iteration, we first consider a Bellman update for the value function of the adversary following \eqref{eq:bellman-recursion-adversary}, and then we compute a Bellman update for the value function of the decision maker following \eqref{eq:pistar:argmax}:
\begin{equation}\label{alg:VI-PI}\tag{VI-PI}
\bm{v}^{0} = \bm{0}, v^{k+1}_{s} = \max_{a \in \A}  \left\{ \; r_{sa} + \lambda \cdot \sum_{i=1}^{r} u^{i}_{sa}\min_{\bm{w}_{i} \in \W^{i}} \; \bm{w}_{i}^{\top}\bm{v}^{k}  \; \right\}, \forall \; s \in \X, \forall \; k \geq 0.
\end{equation}
 We would like to recall that in this paper, (PI) refers to the \textit{Policy Improvement},  so that Algorithm \ref{alg:VI-PI} is not related to policy iteration.
We can now state the main theorem of our paper.
\begin{theorem}\label{th:main}
Let $\PP$ be an $r$-rectangular uncertainty set.
\tb{Under Assumption \ref{ass:tractable-u-set},} Algorithm \ref{alg:VI-PI} returns an $\epsilon$-optimal solution to the policy improvement problem in \tb{$O\left( rSA \log\left(\epsilon^{-1}\right) + \sum_{i=1}^{r} comp\left(\W^{i}\right)  \log^{2}\left(\epsilon^{-1} \right)\right)$} arithmetic operations.
\end{theorem}
\proof{Proof.}
From Equation \eqref{eq:argminargmax} in Theorem \ref{th:duality}, we have the following two coupled Bellman equations:
\[ v^{*}_{s} = \max_{\pi_{s} \in \Delta} \left\{ r_{\pi_{s},s} + \lambda \cdot (\bm{T}_{\pi} \bm{W}^{* \; \top} \bm{v}^{*})_{s} \right\}, \forall \; s \in \X,\]
\[ \beta^{*}_{i} = \min_{\bm{w}_{i} \in \W^{i}} \left\{ \bm{w}_{i}^{\top} ( \bm{r}_{\pi^{*}}+\lambda \cdot \bm{T}_{\pi^{*}} \bm{\beta}^{*}) \right\}, \forall \; i \in [r].\]
From Lemma \ref{lem:Wv-beta} we have $\bm{W}^{* \; \top}\bm{v}^{*} = \bm{\beta}^{*}$, {\color{black} and we can replace $\bm{W}^{* \; \top}\bm{v}^{*}$ in the first equality by $\bm{\beta}^{*}.$ This leads to the following fixed-point equality:}
\begin{equation}
\label{eq:intric-fix-pt}
\begin{aligned}
v^{*}_{s}  = \max_{\pi_{s} \in \Delta} \{ r_{\pi_{s},s} +  \lambda \cdot (\bm{T}_{\pi}(\min_{\bm{w}_{i} \in \W^{i}} \{ \bm{w}_{i}^{\top}\bm{v}^{*} \})_{i \in [r]})_{s} \}, \forall \; s \in \X.
\end{aligned} 
\end{equation}
Since we can choose an optimal robust policy that is deterministic, note that \eqref{eq:intric-fix-pt} is equivalent to
\[
v^{*}_{s} = \max_{a \in \A} \left\{  r_{sa} + \lambda \cdot \sum_{i=1}^{r} u^{i}_{sa}\min_{\bm{w}_{i} \in \W^{i}} \; \bm{w}_{i}^{\top}\bm{v}^{*} \right\}, \forall \; s \in \X.\]
We define the function $F:\R^{S} \rightarrow \R^{S}$ as
\[
 F(\bm{v})_{s} = \max_{\pi_{s} \in \Delta} \left\{ r_{\pi_{s},s} +  \lambda \cdot (\bm{T}_{\pi}(\min_{\bm{w}_{i} \in \W^{i}} \{ \bm{w}_{i}^{\top}\bm{v} \})_{i \in [r]})_{s} \right\}, \forall \; s \in \X.
\]
 The function $F$ is a component-wise non-decreasing contraction, see Appendix \ref{app:lem:F1}. Therefore, its fixed-point is the unique solution to the optimality Equation \eqref{eq:intric-fix-pt}.
To solve the policy improvement problem, it is sufficient to compute the fixed-point of $F$. Following~\cite{Puterman}, we know that the condition $\|\bm{v}^{k+1} - \bm{v}^{k}\|_{\infty} < \epsilon(1-\lambda)(2 \lambda)^{-1}$
is sufficient to conclude that
$\| \bm{v}^{k+1}-\bm{v}^{*} \|_{\infty} \leq \epsilon.$ Therefore, 
Algorithm \ref{alg:VI-PI} returns an $\epsilon$-optimal solution to the policy improvement problem. We now present the analysis of the running time of Algorithm \ref{alg:VI-PI}.

\vspace{2mm}
\noindent \tb{{\textbf{Running time of Algorithm \ref{alg:VI-PI}.}}}
\begin{itemize}
\item \tb{\textit{Number of iterations of Algorithm \ref{alg:VI-PI}.} Following Theorem 6.3.3 in \cite{Puterman}, we can stop Algorithm \ref{alg:VI-PI} as soon as $\|\bm{v}^{k+1} - \bm{v}^{k}\|_{\infty} < \epsilon(1-\lambda)(2 \lambda)^{-1}$ to ensure that $\| \bm{v}^{k+1} - \bm{v}^{*} \|_{\infty} \leq \epsilon$, and this condition will be satisfied after at most 
$O ( \log(\epsilon^{-1}))$ iterations of Algorithm \ref{alg:VI-PI}.}
\item \tb{\textit{Computational cost of each iteration.}
We can start each iteration by solving $\min_{\bm{w}_{i} \in \W^{i}} \bm{w}_{i}^{\top}\bm{v}^{k}$ up to accuracy $\epsilon$, for all $i \in [r]$. This can be done in $O\left(\sum_{i=1}^{r} comp\left(\W^{i}\right) \log(\epsilon^{-1}) \right)$ arithmetic operations. Then for each state $s \in \X$ and action $a \in \A$, we compute 
$r_{sa} + \lambda \cdot \sum_{i=1}^{r} u^{i}_{sa}\min_{\bm{w}_{i} \in \W^{i}} \; \bm{w}_{i}^{\top}\bm{v}^{k}$ in $O(r)$ operations, and then we find, for each $s \in \X$, 
\[\max_{a \in \A}  \{ \; r_{sa} + \lambda \cdot \sum_{i=1}^{r} u^{i}_{sa}\min_{\bm{w}_{i} \in \W^{i}} \; \bm{w}_{i}^{\top}\bm{v}^{k}  \; \}.\]}
\tb{
\noindent
Therefore, the total number of arithmetic operations for each iteration is in $O\left(rSA + \sum_{i=1}^{r} comp(\W^{i}) \log(\epsilon^{-1}) \right) $. }
\end{itemize} 
\noindent
\tb{
Overall, Algorithm \ref{alg:VI-PI} returns an $\epsilon$-optimal policy in $O\left( rSA \log\left(\epsilon^{-1}\right) + \sum_{i=1}^{r} comp\left(\W^{i}\right)  \log^{2}\left(\epsilon^{-1} \right)\right)$ arithmetic operations.}
 \hfill \qed
\endproof
\tb{
\subsection{Complexity for polyhedral uncertainty sets.}
To provide more intuition on the complexity of Algorithm \ref{alg:VI-PE} and Algorithm \ref{alg:VI-PE}, we present here the complexity of each algorithm in the case of {\it polyhedral} uncertainty sets. In particular, in this section we assume that there exists an integer $m \in \N$, such that for each $i \in [r]$ we have \[ \W^{i} = \left\{ \; \bm{w} \in \R^{S} \; | \; \bm{A}_{i}\bm{w} \geq \bm{b}_{i}, \bm{w}_{i} \geq \bm{0} \right\},\] where $\bm{b}_{i}$ are vectors  of size $m \in \N$ and $\bm{A}_{i}$ are matrices in $ \R^{m \times S}$. Polyhedral uncertainty sets recover important uncertainty sets, including sets  based on $\ell_{1}$-norms~\citep{Iyengar,Ho}, $\ell_{\infty}$-norm~\citep{givan1997bounded}, and budget of uncertainty sets~\citep{BS2004}. In this case, we have, for each $i \in [r]$, 
\[ comp(\W^{i}) = (m+S)^{3/2}S^{2},\]
i.e., the decision maker can optimize a linear form on $\W^{i}$ up to precision $\epsilon$ in $O\left((m+S)^{3/2}S^{2} \log(\epsilon^{-1})\right)$ with interior point methods~(see Section 4.6.1 in \cite{BenTal-Nemirovski}). }

\vspace{2mm}
\tb{
\noindent \textbf{Complexity of Algorithm \ref{alg:VI-PE}.} In the case of polyhedral uncertainty sets, Algorithm \ref{alg:VI-PE} returns an $\epsilon$-optimal solution to the policy evaluation problem in $O\left(rSA + r \cdot (m+S)^{3/2} S^{2} \cdot \log^{2}\left(\epsilon^{-1}\right)\right)$ arithmetic operations.}

\vspace{2mm}
\tb{
\noindent \textbf{Complexity of Algorithm \ref{alg:VI-PI}.} In the case of polyhedral uncertainty sets, Algorithm \ref{alg:VI-PI} returns an $\epsilon$-optimal solution to the policy improvement problem in $O\left(rSA \log(\epsilon^{-1}) + r \cdot (m+S)^{3/2} S^{2} \cdot \log^{2}\left(\epsilon^{-1}\right)\right)$ arithmetic operations.}

\vspace{2mm}
\noindent \textbf{Comparison with robust value iteration for $(s,a)$-rectangular uncertainty sets.} 
 Note the close relation between Algorithm \ref{alg:VI-PI} for robust MDPs with $r$-rectangular uncertainty set, and the classical Value Iteration \ref{alg:VI-sarec} for robust MDPs with $(s,a)$-rectangular uncertainty sets. If $\PP = \times_{(s,a) \in \X \times \A} \PP_{sa}$,  Algorithm \ref{alg:VI-sarec} runs as
\begin{equation}\label{alg:VI-sarec}\tag{VI-(s,a)}
\bm{v}^{0}=\bm{0},v^{k+1}_{s} = \max_{a \in \A} \left\lbrace r_{sa} + \lambda \min_{\bm{P}_{sa} \in \PP_{sa}} \bm{P}_{sa}^{\top}\bm{v}^{k}, \right\rbrace,  \forall \; s \in \X,  \forall \; k \geq 0.
\end{equation}
To compare with the case of polyhedral $r$-rectangular uncertainty sets, we also assume here that the sets $\PP_{sa}$ are polyhedral with $m$ constraints, for all $(s,a) \in \X \times \A$.

In this case, Algorithm \ref{alg:VI-sarec} returns an $\epsilon$-optimal policy to the robust MDP problem with $(s,a)$-rectangular uncertainty sets after \tb{$O\left(SA \log(\epsilon^{-1})+ (m+S)^{3/2}S^{3}\log^{2}(\epsilon^{-1})  \right)$}. The analysis is similar to the analysis of the complexity of Algorithm \ref{alg:VI-PI}.  This complexity is to be compared with the complexity of Algorithm \ref{alg:VI-PI} given in Theorem \ref{th:main}. Typically, when $A$ is in the same (or smaller) order of magnitude as $S$,  the leading term in the complexity of Algorithm \ref{alg:VI-PI} is \tb{$O\left(r(m+S)^{3/2}S^{2}\log^{2}(\epsilon^{-1})\right)$}, and the leading term in the complexity of Algorithm \ref{alg:VI-sarec} is \tb{$O \left( (m+S)^{3/2}S^{3}\log^{2}(\epsilon^{-1}) \right)$}.  Therefore, if $r$ is smaller than $S$, then Algorithm \ref{alg:VI-PI} enjoys better worst-case guarantees than Algorithm \ref{alg:VI-sarec}.

Additionally,  remember that $r$-rectangularity generalizes $(s,a)$-rectangularity.  In particular, we can convert an $(s,a)$-rectangular uncertainty to an $r$-rectangular uncertainty set, leading to a choice of $r=SA$ (see Proposition \ref{prop:relationship-s-r-sa-rec}).  Such a reformulation would lead to a complexity of \tb{$O\left( S^{2}  A^{2}  \log\left(\epsilon^{-1}\right)+ (m+S)^{3/2}S^{3} A\log^{2}\left(\epsilon^{-1}\right)  \right)$} for Algorithm \ref{alg:VI-PI}, which is worse than the complexity of Algorithm \ref{alg:VI-sarec}.  Note that the complexity of Algorithm \ref{alg:VI-PI} could be improved in this case, since the coefficients $u_{sa}^{i}$ are $0$ or $1$ in this reformulation. Exploiting this in our implementation, Algorithm \ref{alg:VI-PI} would then be exactly equivalent to Algorithm \ref{alg:VI-sarec}. Therefore, the reformulation of $(s,a)$-rectangular uncertainty set as $r$-rectangular uncertainty set provides interesting insights on the modeling power of $r$-rectangular uncertainty sets but does not lead to more efficient algorithms for $(s,a)$-rectangular uncertainty sets.
\subsection{Lifted MDP.}
We would like to note that we can \textit{lift} a robust MDP instance with an $r$-rectangular uncertainty set to a larger MDP instance with an $(s,a)$-rectangular uncertainty set. This lifting is obtained by augmenting the state space and it preserves the optimal policy. \tb{Note that the lifted MDP is $(s,a)$-rectangular only because the state space is larger than in the original $r$-rectangular MDP. Therefore, the construction of the lifted MDP does not contradict the third statement of Proposition \ref{prop:relationship-s-r-sa-rec}.} We present the details below.

\vspace{2mm}
\noindent \textbf{State augmentation.}
Let $\M = \left( \X,\A,\PP,\bm{r},\lambda \right)$ be a robust MDP instance and assume that $\PP$ is $r$-rectangular, with the set of factors being $\W$ and the coefficients being $\left(u_{sa}^{i}\right)_{i \in [r]}$ for each pair $(s,a) \in \X \times \A.$
We consider the following {\it lifted MDP} $\hat{\M} =  \left( \hat{\X},\hat{\A},\hat{\PP},\hat{\bm{r}},\hat{\lambda} \right)$ with an $(s,a)$-rectangular uncertainty set $\hat{\PP}$ as follows.
\begin{itemize}
\item The set of states $\hat{\X}$ is $\hat{\X} = \X \bigcup [r]$, with cardinality $S+r$. The set $\hat{\X}$ is called the \textit{augmented state space}.
\item The set of actions available at state $\hat{s} \in \hat{\X}$ is $\A$ if $\hat{s} \in \X$ and $\{ 1 \}$ if $\hat{s} \in [r]$ (i.e., there is no action to choose if $\hat{s} \in [r]$). The cardinality of $\hat{\A}$ is $A+1$.
Therefore, the set of policies in $\hat{\M}$ is \[ \hat{\Pi} = \left\lbrace \phi(\pi)  \; \bigg| \; \pi \in \Pi \right\rbrace,\]
for $\phi: \R^{S \times A} \rightarrow \R^{(S+r) \times (A+1)}$ defined as
\[\phi(\pi)=\begin{pmatrix}
\pi & \bm{0}_{S}\\
\bm{0}_{r \times A}  & \bm{e}_{r}
\end{pmatrix} \in \R^{(S+r) \times (A+1)}\]  with $\bm{e}_{r}=(1,...,1) \in \R^{r}$. Note that $\phi$ is a one-to-one mapping between the policies in $\hat{\M}$ and the policies in $\M$.
\item For $(s,a) \in \X \times \A$, we have $\hat{r}_{sa} = r_{sa}$, and for $i \in [r]$ we have $r_{i,1}=0$.
\item For $(s,a) \in \X \times \A$,  we have $\hat{\bm{P}}_{sa}  = \left\lbrace \left( \bm{0}_{S}; \left( u_{sa}^{i} \right)_{i \in [r]} \right) \right\rbrace$.  For $i \in [r]$, $\hat{\bm{P}}_{i,1} = \left\lbrace  \left( \bm{w}_{i} ; \bm{0}_{r} \right) \; | \; \bm{w}_{i} \in \W^{i} \right\rbrace.$
\item The discount factor is $\hat{\lambda} = \sqrt{\lambda}$.
\end{itemize}
 \tb{A comparison of the dynamics for the $r$-rectangular MDP $\M$ and the lifted MDP $\hat{\M}$ is presented in Figure \ref{fig:comparison-factor-lifted-mdp}. Since in the $r$-rectangular MDP we have
\begin{equation}\label{eq:link-psa-u-w}
P_{sas'} = \sum_{i=1}^{r} u_{sa}^{i}w_{i,s'}, \forall \; s,s' \in \X, a \in \A,
\end{equation}
and since
\[
\sum_{i=1}^{r} u_{sa}^{i} =1, \sum_{s' \in \X} w_{i,s'} =1, \forall \; s \in \X, a \in \A, i \in [r],\]
we can interpret Equality \eqref{eq:link-psa-u-w} as an intermediate transition from the pair $(s,a)$ to a factor $i \in [r]$ (with probability $u_{sa}^{i}$) followed by a transition to the next state $s'$ (with probability $w_{i,s'}$). In short, the lifted MDP $\hat{\M}$ integrates the set $[r]$ to the set of states $\X$. The transitions from $\X$ to $[r]$ happen during even periods $\{ 2 t| t \in \N\}$ and the transitions from $[r]$ to $\X$ happen during odd periods $\{2t+1|t \in \N\}$, as illustrated in Figure \ref{fig:lifted_mdp}.}

\tb{Note that in the lifted MDP, the transitions from $s \in \X$ to $i \in [r]$ are not subject to parameter uncertainty. Only the transitions from $i$ to $s'$ (which happen with probability $w_{i,s'}$) are uncertain. If the robust MDP $\M$ is $r$-rectangular, then the factors $\bm{w}_{1},...,\bm{w}_{r}$ associated with $i=1,...,r$ can vary independently. Since in the lifted MDP $\hat{\M}$ the set $[r]$ is a subset of the set of states, and since no action is chosen in the state $i=1,...,r$, the $r$-rectangularity of the robust MDP $\M$ becomes an $(s,a)$-rectangularity assumption in the lifted MDP $\hat{\M}$.}
\tb{
\begin{figure}[h]
\begin{subfigure}{0.4\textwidth}
  \includegraphics[width=0.9\linewidth]{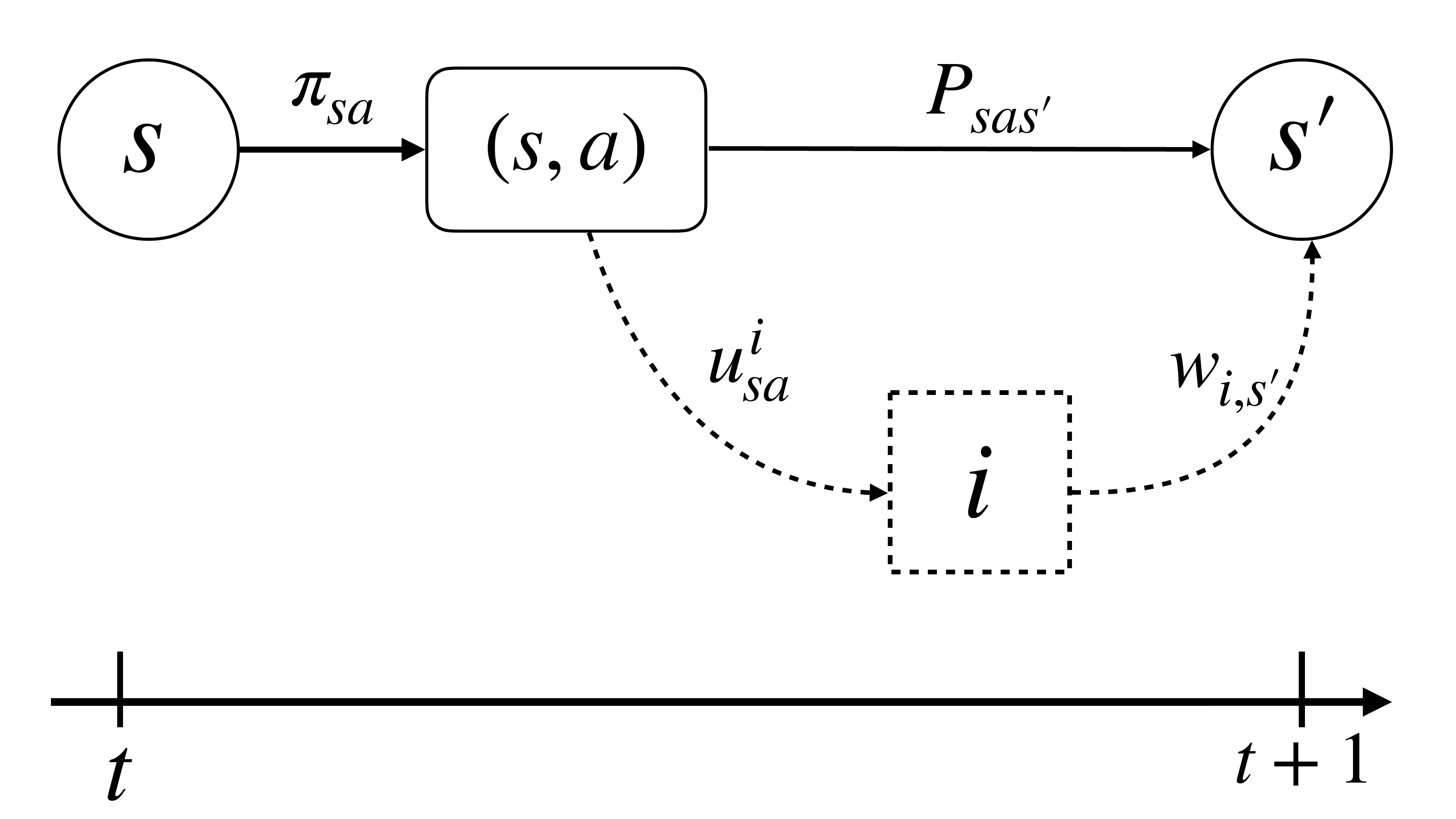}
\caption{Dynamics from $s$ to $s'$ in the robust MDP $\M$.}
\label{fig:factor_mdp}
\end{subfigure}
\begin{subfigure}{0.55\textwidth}
  \includegraphics[width=0.9\linewidth]{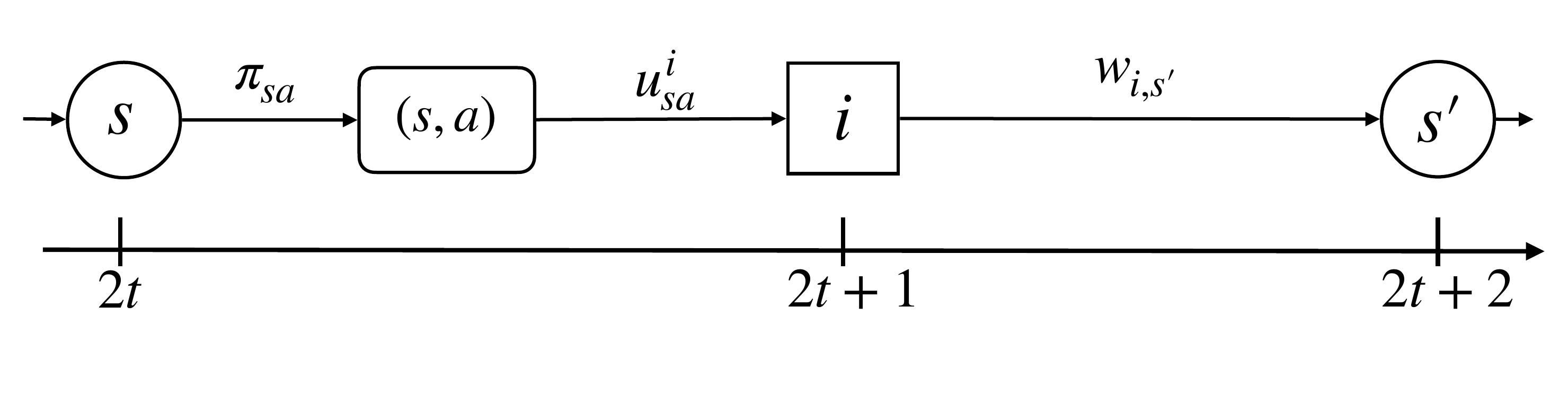}
\caption{Dynamics from $s$ to $i$ to $s'$ in the lifted MDP $\hat{\M}$.}
\label{fig:lifted_mdp}
\end{subfigure}
\caption{\tb{Comparison of the dynamics for the robust MDP with factor matrix uncertainty set (Figure \ref{fig:factor_mdp}) and for the corresponding lifted MDP (Figure \ref{fig:lifted_mdp}). The states $s,s' \in \X$ are represented with circles. The factor $i \in [r]$ is represented with a square. The transition probabilities are represented above the transition arcs. The dashed transition arcs in Figure \ref{fig:factor_mdp} make us of the decomposition $P_{sas'} = \sum_{i=1}^{r}u^{i}_{sa}w_{i,s'}$ from the definition of the factor matrix uncertainty set in \eqref{eq:fmus}. In the lifted MDP of Figure \ref{fig:lifted_mdp}, the states $s,s' \in \X$ are attained at periods $2t$ and $2t+2$, and the factor $i \in [r]$ is attained at period $2t+1$.}}
\label{fig:comparison-factor-lifted-mdp}
\end{figure}}

\vspace{2mm}
\noindent \tb{ \textbf{Value function and optimal policies.}
Let $\hat{\pi} \in \hat{\Pi}$ be a policy for the decision maker in the lifted MDP $\hat{\M}$ and let $\hat{\bm{v}} \in \R^{S+r}$ be its value function. By definition of $\hat{\Pi}$, we have $\hat{\pi} = \phi(\pi)$ for some $\pi \in \Pi$. For $s \in \X$, we have
\begin{align}
\hat{v}_{s} & = \E^{\hat{\pi}, \hat{\bm{P}}} \left[ \sum_{t=0}^{\infty} \hat{\lambda}^{t}\hat{r}_{\hat{s}_{t}\hat{a}_{t}} \; \bigg| \; \hat{s}_0 = s \right] \label{eq:value-vec-1} \\
& = r_{\pi,s} + \sqrt{\lambda} \sum_{a \in \A} \sum_{i=1}^{r} \pi_{sa} u_{sa}^{i} \cdot \E^{\pi, \bm{P}} \left[ \sum_{t=0}^{\infty} \hat{\lambda}^{t}\hat{r}_{\hat{s}_{t}\hat{a}_{t}} \; \bigg| \; \hat{s}_0 = i \right]\label{eq:value-vec-2} \\
& = r_{\pi,s} + \sqrt{\lambda} \sum_{a \in \A} \sum_{i=1}^{r} \pi_{sa} u_{sa}^{i} \sum_{s' \in \X} w_{i,s'} \sqrt{\lambda} \cdot \E^{\pi, \bm{P}} \left[ \sum_{t=0}^{\infty} \hat{\lambda}^{t}\hat{r}_{\hat{s}_{t}\hat{a}_{t}} \; \bigg| \; \hat{s}_0 = s' \right]\label{eq:value-vec-3} \\
& = r_{\pi,s} + \lambda \sum_{a \in \A} \sum_{i=1}^{r} \pi_{sa} u_{sa}^{i} \sum_{s' \in \X} w_{i,s'}  \hat{v}_{s'}\label{eq:value-vec-4}\\
& = r_{\pi,s} + \lambda \sum_{a \in \A} \pi_{sa} \bm{P}^{\top}_{sa}  \left(\hat{v}_{s'}\right)_{s' \in \X}\label{eq:value-vec-5}
\end{align}
where \eqref{eq:value-vec-1} follows from the Bellman equation in $\hat{\M}$, and \eqref{eq:value-vec-2} and \eqref{eq:value-vec-3} follow from the definition of the rewards in $\hat{\M}$. Equality \eqref{eq:value-vec-5} follows from $P_{sas'} = \sum_{i=1}^{r} u_{sa}^{i}w_{i,s'}$ in a factor MDP. Note that \eqref{eq:value-vec-5} is exactly the Bellman equation for $\M$. 
Therefore, the value function in $\hat{\M}$ at state $s \in \X$ is equal to the value function in $\M$ at the same state $s \in \X$. Let us write $\hat{R}(\hat{\pi},\hat{\bm{P}})$ the reward of the decision maker in the lifted MDP $\hat{\M}$. Equality \eqref{eq:value-vec-5} shows that if $\hat{\pi} = \phi(\pi)$, then 
\[\hat{R}(\hat{\pi},\hat{\bm{P}}) = R(\pi,\bm{P}),\]
i.e., the lifted MDP $\hat{\M}$ preserves optimality of policies from the robust MDP $\M$.}

\vspace{2mm}
\noindent \textbf{Comparison of the reformulation and our approach.}
Some of our results on the structure of optimal policies follow from this lifting and the properties of $(s,a)$-rectangular uncertainty sets. In particular, Lemma \ref{lem:determ-poli} (optimality of stationary deterministic policy) follows from Theorem 4 in \cite{Iyengar}. Additionally,  the duality result \eqref{eq:strong-min-max-duality} in Theorem \eqref{th:duality} follows from Theorem 3 in \cite{Nilim},  although our stronger equilibrium result \eqref{eq:argminargmax} is new.
The correctness of Algorithm \ref{alg:VI-PE} follows from Theorem 2 in \cite{Kuhn} (along with our Lemma \ref{lem:PEfinal}), while the correctness of  Algorithm \ref{alg:VI-PI}  follows from Theorem 3 in \cite{Nilim}.

Even though this state augmentation simplifies some of the results presented in Section \ref{sec:duality} and Section \ref{sec:PI},  there are some advantages in the more direct approach that we present in this paper. 
First, working directly with the MDP $\M$ instead of $\hat{\M}$ simplifies the formulation of the robust value iteration algorithm, since the state space in $\M$ is smaller than in $\hat{\M}$.
Additionally,  the intuition that the adversary from Section \ref{sec:PE} is solving an MDP is important and leads to our new Nash equilibrium result \eqref{eq:argminargmax}. 
We then use this Nash equilibrium result to write the two coupled Bellman equations for the decision maker and the adversary.  This leads to the important intuition that the robust MDP problems with rectangular uncertainty sets is a two-player game {\it where each player plays a separate MDP instance}. This Nash equilibrium result is also crucial to prove the robust maximum principle and our robust Blackwell optimality result (see next section).  
\section{Properties of optimal robust policies.}\label{sec:theory} For $r$-rectangular uncertainty sets, we are able to prove important structural properties of optimal robust policies.

\subsection{Robust maximum principle.}
 For a classical MDP, the value function of the optimal policy is component-wise higher than the value function of any other policy. This is known as the \textit{maximum principle} (\cite{handbookMDP}, Section 2). {\color{black} Therefore, an optimal nominal policy attains the highest possible nominal value, \textit{starting from each state}: this is a special case of Pontryagin's maximum principle, which is at the core of the theory of optimal control \citep{ross2009primer}.}
We prove a \textit{robust maximum principle} for robust MDPs. { \color{black} The robust maximum principle has been leveraged to prove structural properties of optimal solutions of MDPs in a healthcare setting \citep{grand2020robust},  since it provides useful inequalities to compare the optimal nominal and the optimal robust policies. In particular, it shows that the optimal robust policy achieves the highest possible worst-case value, \textit{starting from any state}.}   
We write  $\bm{v^{\pi,\bm{W}}}$ the value function of the decision maker when s/he chooses policy $\pi$ and the adversary chooses factor matrix $\bm{W}.$

\begin{proposition}\label{prop:robust-max-principle} \tb{Under Assumption \ref{ass:tractable-u-set},} let $\PP$ be an $r$-rectangular uncertainty set.
\begin{enumerate}
\item Let $\pi$ be a policy and $\bm{W}^{1} \in \arg \min_{\bm{W} \in \W} \; R(\pi,\bm{W}).$
Then
\[v^{\pi,\bm{W}^{1}}_{s} \leq v^{\pi,\bm{W}^{0}}_{s},\;\forall \; \bm{W}^{0} \in \W, \forall \; s \in \X.\]
\item Let $(\pi^{*},\bm{W}^{*}) \in \arg \; \max_{\pi \in \Pi} \; \min_{\bm{W} \in \W} \; R(\pi,\bm{W}).$ Then 
\tb{
\[\forall \; \pi \in \Pi, \forall \; \bm{W}^{1} \in \arg \min_{\bm{W} \in \W} \; R(\pi,\bm{W}), v^{\pi,\bm{W}^{1}}_{s} \leq v^{\pi^{*},\bm{W}^{*}}_{s}, \forall \; s \in \X.\]}
\end{enumerate}
\end{proposition}
We present a detailed proof in Appendix \ref{app:robust-max-principle}. The proof relies on Lemma \ref{lem:Wv-beta} and the strong duality of Theorem \ref{th:duality}.  In particular, because Theorem \ref{th:duality} does not necessarily hold for some $s$-rectangular uncertainty sets, our proof cannot be adapted to an $s$-rectangular uncertainty set.

\subsection{Robust Blackwell optimality.}
We first review the notion of Blackwell optimality for nominal MDPs.
In the classical MDP literature, given a fixed known transition kernel $\bm{P}$, a policy $\pi$ is said to be \textit{Blackwell optimal} if it is optimal for all discount factors close enough to 1 (Section 10.1.2 in \cite{Puterman}).  Blackwell optimal policies are optimal in the long-run average reward MDP, obtained by changing the reward criterion from the discounted cumulated reward \eqref{eq-expreward} to the long-run average reward: $\lim_{T \rightarrow + \infty} \dfrac{1}{T} \E^{\pi, \bm{P}} \left[ \sum_{t=0}^{T} r_{s_{t}a_{t}} \; \bigg| \; s_{0} \sim \bm{p}_{0} \right]$. However,  the long-run average criterion focuses on the  \textit{steady-state} performance of the policy, and does not reflect any reward gathered in finite time, i.e., it does not take into account the \textit{transient} performance of the policy.  As a consequence, there may be multiple optimal policies for the long-run average reward criterion. Hence, there may be a preference for a more selective criterion, while still accounting for the long-run average reward. 
For instance, one may ask for an $n$-discount optimal policy $\pi_{n}$, which for some integer $n \geq 1$ satisfies
\begin{equation}\label{eq:n-discount-optimal}
\lim_{\lambda \rightarrow 1} \left(1-\lambda\right)^{-n} \left(\bm{v}^{\pi_{n}}_{\lambda} - \bm{v}^{\pi}_{\lambda} \right) \geq 0, \forall \; \pi \in \Pi,
\end{equation}
where $\bm{v}^{\pi}_{\lambda}$ is the value function of policy $\pi \in \Pi$ for the discount factor $\lambda$. Note that the selectivity of the optimality criterion \eqref{eq:n-discount-optimal} increases with $n$ (see \cite{Puterman}, Section 10.1.1 for more details and relations to \textit{bias-optimality}). 
A Blackwell optimal policy is $n$-discount optimal, \textit{for any integer} $n \geq 1$. Therefore,  it is still optimal for the long-run average criterion, but it is much more sensitive to the reward obtained after any finite numbers of periods.

We extend the notion of Blackwell optimality for robust MDP where the uncertainty set is $r$-rectangular.  
To the best of our knowledge, there is no results on robust MDPs with the long-run average reward criterion. The closest problems are two-player perfect information zero-sum stochastic games with finite number of actions, for which the existence of Blackwell optimal policies is proved in \cite{gaubert1998non}. As a special case,  and without explicitly making the connection to robust MDP, the authors in \cite{akian2019operator} obtain the existence of Blackwell optimal policies in an instance of an $s$-rectangular robust MDP with long-run average reward (\cite{akian2019operator}, Theorem 8), where the rewards are Kullback-Leibler divergences between the policy $\pi_{s}$ chosen at a state $s$ and some predefined weight vectors. 
We do not make any assumption on the structure of the rewards $r_{sa}$, and we assume that $\PP$ is an $r$-rectangular uncertainty set. \tb{However, we need to assume that the sets $\W^{1},...,\W^{r}$ have finitely many extreme points. For instance, this is the case when $\W^{1},...,\W^{r}$ are polyhedral.  In particular, we have the following proposition.  We present a detailed proof of Proposition \ref{prop:black} in the Appendix \ref{app:black}. }

\begin{proposition}\label{prop:black}
\tb{Let $\PP$ be an $r$-rectangular uncertainty set. Assume that the sets $\W^{1}, ..., \W^{r}$ have finitely many extreme points.
Then there exists a stationary deterministic policy $\pi^{*}$, a factor matrix $\bm{W}^{*}$, and a discount factor $\lambda_{0} \in (0,1),$ such that for all $\lambda \in (\lambda_{0},1)$,
the pair $(\pi^{*},\bm{W}^{*})$ remains an optimal solution to the robust MDP problem \eqref{eq:robust-mdp-problem}.}
\end{proposition}
We extend Proposition \ref{prop:black} to the interval $[0,1]$ in the next proposition.  The proof of Proposition \ref{prop:blackwell-01} is in Appendix \ref{app:black-01}.
\tb{
\begin{proposition}\label{prop:blackwell-01}
Let $\PP$ be an $r$-rectangular uncertainty set. Assume that the sets $\W^{1}, ..., \W^{r}$ have finitely many extreme points. Then there exists an integer $p \in \N$, there exists some scalars $\lambda_{0}=0 < \lambda_{1} < ... < \lambda_{p}=1$ such that  for all $j \in \{0,...,p-1\}$,  the same pair of stationary deterministic policy and factor matrix $(\pi_{j},\bm{W}_{j})$ is an optimal solution to the robust MDP problem \eqref{eq:robust-mdp-problem} for all $\lambda \in (\lambda_{j},\lambda_{j+1})$.
\end{proposition}
}
Note that for nominal MDPs,  the same proposition as Proposition \ref{prop:blackwell-01} holds and it is possible to compute the breakpoints $\lambda_{0}, ..., \lambda_{p}$ using sensitivity analysis and the simplex algorithm \citep{hordijk1985sensitivity}. Crucially, this algorithm is based on the linear programming reformulation of nominal MDPs. However, no such reformulation is known for robust MDPs.  Our proof of Proposition \ref{prop:blackwell-01} is not constructive,  and it is an interesting open question to find an algorithm that can efficiently compute the breakpoints $\lambda_{0}, ..., \lambda_{p}$ in the case of robust MDPs.
\section{Numerical experiments.}\label{sec:simu}
In this section we study the numerical performance of our optimal robust policies. We show that an optimal nominal policy may be sensitive to small parameter deviations, highlighting the need for designing policies that are robust to parameter uncertainty. We also emphasize the role of nonnegative matrix factorization in constructing factor matrix uncertainty sets. We compare the empirical performances of robust policies for $r$-rectangular and $s$-rectangular uncertainty sets and show empirically that $r$-rectangular uncertainty sets may lead to robust policies that are less pessimistic than for $s$-rectangular uncertainty sets.
\subsection{Machine replacement problem.}
We start with the following machine replacement MDP.

\vspace{2mm}
\noindent \textbf{Problem setup.}

We consider a machine maintenance model introduced in~\cite{Delage} and \cite{Kuhn}. There are 10 states, $\{1,..., 8\} \bigcup \{R1,R2 \}$. The states $1$ to $8$ model the states of deterioration of a machine. There is a reward of $20$ in states $1$ to $7$ while the reward is $0$ in state $8$. There are two repair states $R1$ and $R2$. The state $R1$ is a normal repair and has reward of $18$, the state $R2$ is a long repair and has reward $10$.  The discount factor is $\lambda = 0.8$. The initial distribution is uniform across all states. There are two actions, \textit{wait} and \textit{repair}, and the goal is to maximize the infinite horizon expected reward. We present the nominal transition rate $\bm{P}^{\sf nom}$ in Figure \ref{fig:Ex_Kuhn}. We assume that we know the nominal kernel $\bm{P}^{\sf nom}$ as well as an upper bound $\tau >0$ on the maximum deviation from any component of $\bm{P}^{\sf nom}$. 
We compute the optimal nominal policy and show its sensitivity to transition probabilities.
We then compare the performances of the robust policies associated with the $r$-rectangular and the $s$-rectangular uncertainty sets.
\begin{figure}[h]
\begin{subfigure}{0.5\textwidth}
  \includegraphics[width=0.9\linewidth]{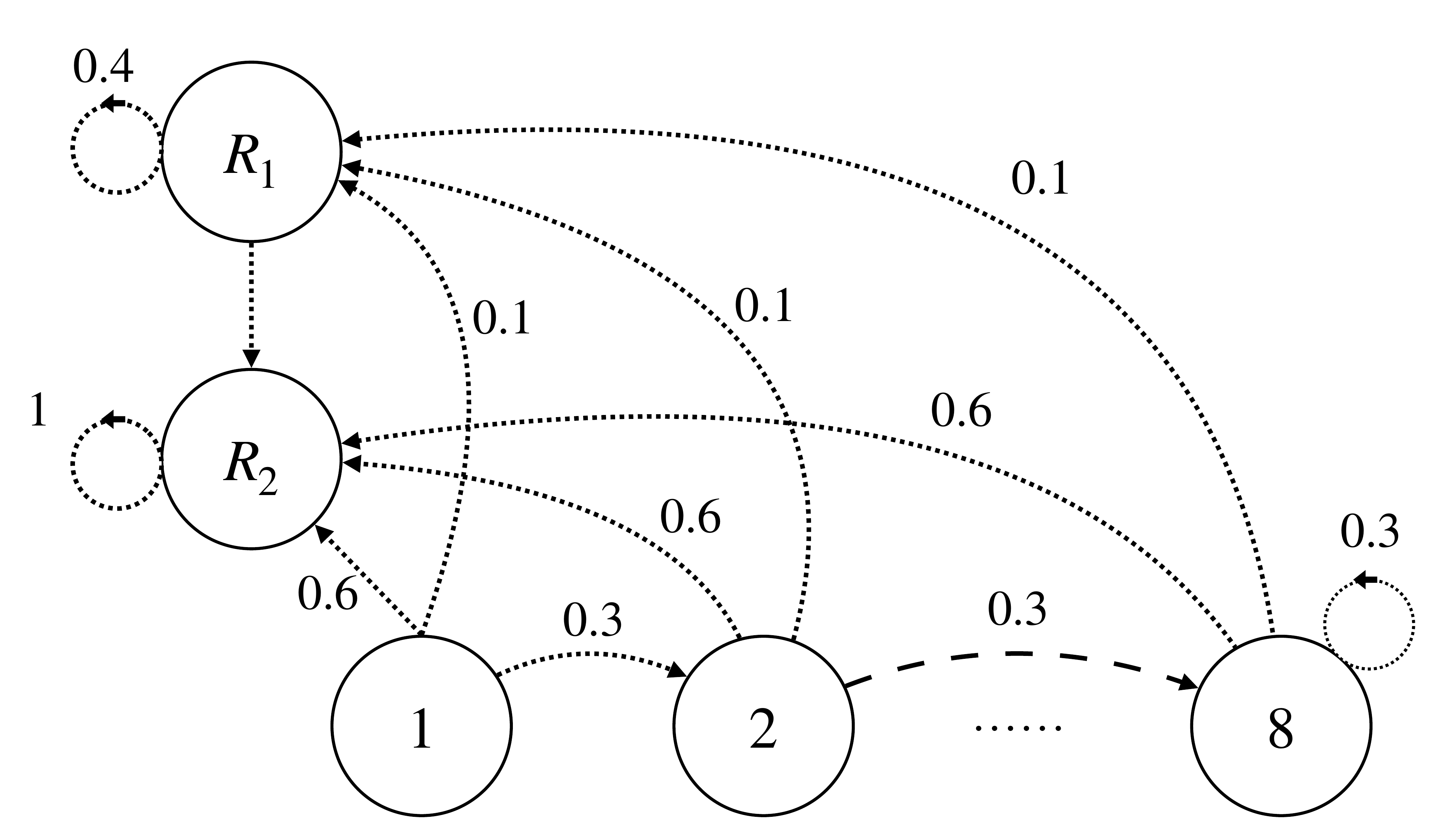}
\caption{Transition probabilities for \\ action \textit{repair}.}
\label{fig:subim1_K}
\end{subfigure}
\begin{subfigure}{0.5\textwidth}
  \includegraphics[width=0.9\linewidth]{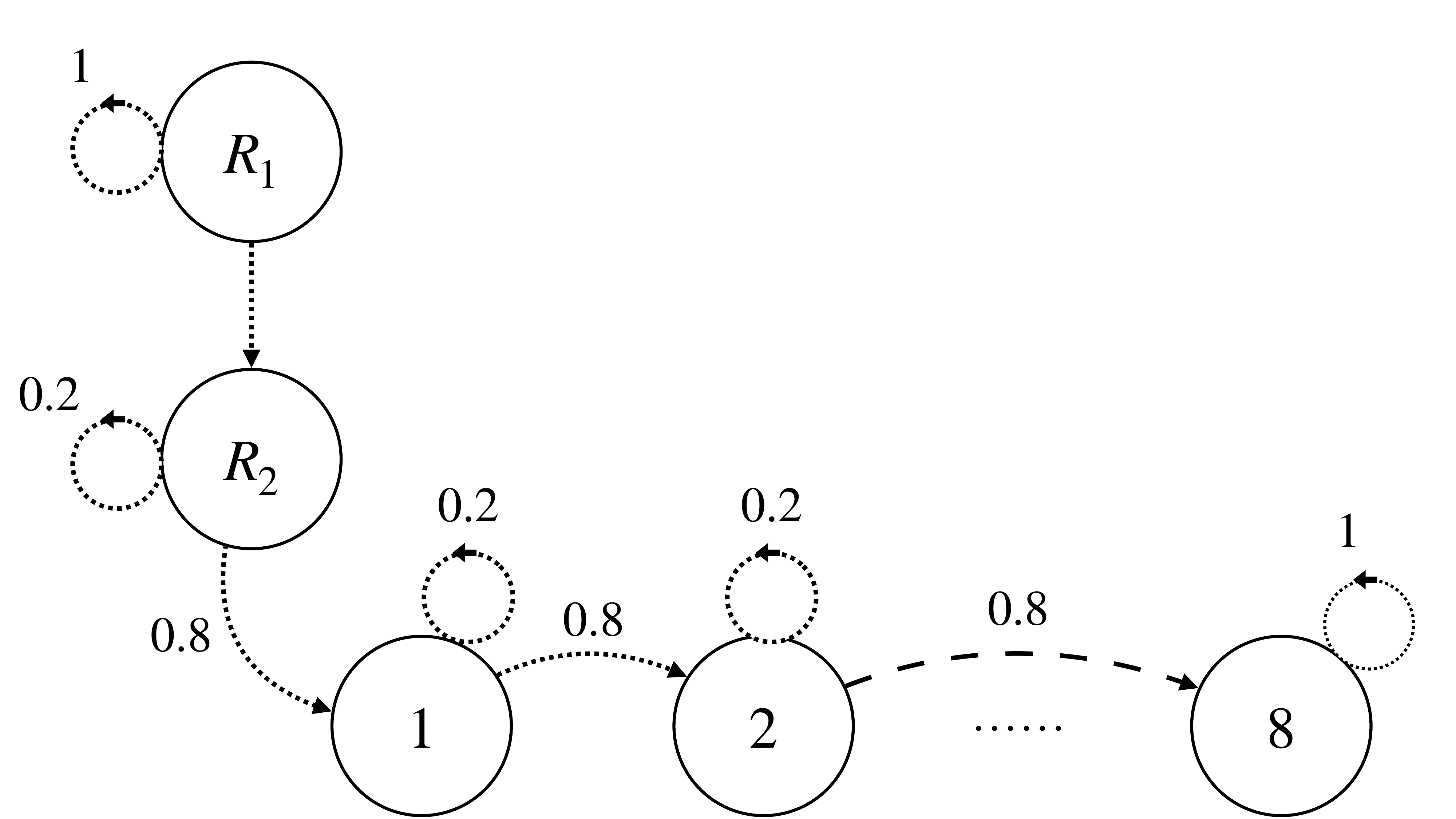}
\caption{Transition probabilities for \\ action \textit{wait}.}
\label{fig:subim2_K}
\end{subfigure}
\caption{Transition probabilities for the machine replacement example. There is a reward of $18$ in state $R1$, of $10$ in state $R2$ and of $0$ in state $8$. All others states have a reward of $20$.  }
\label{fig:Ex_Kuhn}
\end{figure}

\noindent \textbf{Construction of the $r$-rectangular uncertainty set.}
We compute the factor matrix $\bm{W}$ and the coefficients matrices $\bm{u}_{1},...,\bm{u}_{S}.$ Let $\bm{P}^{\sf nom \; \top}$ be the matrix whose columns correspond to the transitions $\bm{P}^{\sf nom}_{sa}$ for every state-action pair $(s,a)$.
The decomposition $\bm{P}^{\sf nom}_{sa} = \sum_{i=1}^{r} u^{i}_{sa}\bm{w}_{i}$ is equivalent to $(\bm{P}_{sa})_{a \in \A} = 
\bm{Wu}_{s},$
 for some factor matrix $\bm{W}$ in $\R_{+}^{S \times r}$ and some coefficients matrices $\bm{u}_{1}, ...,\bm{u}_{S}$ in $\R_{+}^{r \times A}$. Therefore, to estimate $\bm{W}^{\sf nom}, \bm{u}^{\sf nom}$, we solve the following Nonnegative Matrix Factorization (NMF) program: for
 $H_{1} =\{ \bm{W} \; | \; \bm{W} \in \R_{+}^{S \times r}, \bm{W}^{\top}\bm{e}_{S}=\bm{e}_{r} \; \}, H_{2}  = \underset{s \in \X}{\times} \; H_{3} \text{ where } H_{3} = \{ \bm{u} \; | \; \bm{u} \in \R_{+}^{r \times A}, \bm{u}^{\top}\bm{e}_{r}=\bm{e}_{A} \; \}$, we solve
 \begin{align}\label{min:NMF}
 \min \; & \dfrac{1}{2} \| \bm{\tilde{P}}^{\sf nom \; \top} - \bm{Wu} \|_{2}^{2} \\
 & \bm{W} \in H_{1}, \bm{u} \in H_{2}. \label{eq:NMF:constr}
 \end{align}
{\color{black}
Note that a related method has been proposed in \cite{duan2018state}, where the authors propose a state-aggregation mechanism relying on the spectral decomposition of the transition matrix.}

For this example we consider a rank of $r=12$ and we compute a local optimal solution $(\bm{W}^{\sf nom},\bm{u}^{\sf nom})$ of \eqref{min:NMF} by adapting a classical algorithm for NMF~\citep{Wotao}. {\color{black} The choice of $r=12$ enables us to obtain a reasonably good NMF approximation of the transition kernel.  In particular, our solution $(\bm{W}^{\sf nom},\bm{u}^{\sf nom})$ achieves the following errors: if we write \[\bm{M}_{err}=\bm{\tilde{P}}^{\sf nom \; \top} - \bm{W}^{\sf nom} \bm{u}^{\sf nom},\] then
$\| \bm{M}_{err} \|_{2}=7.6\cdot 10^{-4},\| \bm{M}_{err} \|_{1}=2.6 \cdot 10^{-3},\| \bm{M}_{err} \|_{\infty}=2.5 \cdot 10^{-4}.$  Additionally, since the number of states here is $S=10$ and the number of actions is $A=2$, the number of transition vectors $\bm{P}_{sa}$ is $SA=20$. Therefore, the choice $r=12$ represents a rank reduction of $8$.  Smaller choices of $r$ lead to a worst-performance of the NMF solution as an approximation of the nominal transition kernel.}
We use $(\bm{w}_{1}^{\sf nom},...,\bm{w}_{r}^{\sf nom})$ as the nominal factor vectors and we find the coefficients $(\bm{u}_{1},...,\bm{u}_{S})$ as the blocks of the matrix $\bm{u}^{\sf nom}$.

\noindent \textbf{Model of uncertainty.}
We consider the following \textit{budget of uncertainty} set introduced in~\cite{BS2004}: $\PP^{\sf (r)}  = \left\lbrace \left( \sum_{i=1}^{r} u^{i}_{sa}w_{i,s'} \right)_{sas'}  \bigg| \;    \bm{W} = (\bm{w}_{1},..., \bm{w}_{r}) \in \W \right\rbrace$, where $\W  = \W^{1} \times ... \times \W^{r}$ and
\begin{equation*}
\begin{aligned}
\W^{i} & =\{\bm{w}_{i}=\bm{w}_{i}^{\sf nom}+ \bm{\delta} \; | \;  \bm{\delta} \in \R^{S}, \; \| \bm{\delta} \|_{1} \leq \sqrt{S}\cdot \tau,\| \bm{\delta} \|_{\infty} \leq \tau, \; \bm{e}_{S}^{\top}\bm{w}_{i}=1, \; \bm{w}_{i} \geq \bm{0} \}, \; i \in [r].
\end{aligned}
\end{equation*}
The deviations on each component of the factors $\bm{w}_{1},...,\bm{w}_{r}$ are constrained to be smaller than $\tau.$ Moreover, for each factor vector the deviations on each component are independent and the total deviation is constrained to be smaller than $c \sqrt{S}\cdot \tau.$ We refer to \cite{BS2004}, Theorem 3,  and \cite{bertsimas2006robust}, Section 2.1, for details about the choice of the budget of uncertainty,  and the results relating the value of the deviation in $\| \cdot \|_{1}$,  the probability of constraint violation and the degree of conservatism of optimal robust policies.

\noindent \textbf{Construction of the $s$-rectangular uncertainty set.}
We consider the following budget of uncertainty set where the matrices $(\bm{P}_{sa})_{a \in \A} \in \R^{A \times S}$ are not related across different states: $\PP^{\sf (s)}  = \underset{s \in \X}{\times} \; \PP^{\sf (s)}_{s}$, where
\begin{equation*}
\begin{aligned}
\PP^{\sf (s)}_{s}& =\{\bm{P}_{s}=\bm{P}_{s}^{\sf nom} + \bm{\Delta} \; | \;  \bm{\Delta} \in \R^{A \times S}, \; \| \bm{\Delta} \|_{1} \leq \sqrt{S \cdot A} \cdot \tau,\| \bm{\Delta} \|_{\infty} \leq \tau, \; \bm{P}_{s}\bm{e}_{S}=\bm{e}_{A}, \; \bm{P}_{s} \geq \bm{0} \}, s \in \X.
\end{aligned}
\end{equation*}
The maximum deviation from each component $P^{\sf nom}_{sas'}$ is $\tau$. For the same reason as for the $r$-rectangular uncertainty set of the previous section, the total deviation from a given matrix $\bm{P}_{s}^{\sf nom}$ is $ \sqrt{S \cdot A} \cdot \tau$.

\noindent \textbf{Empirical results.}
We compute the optimal nominal policy $\pi^{\sf nom}$ using value iteration (\cite{Puterman}, Chapter 6.3) and we normalize the reward so that $R(\pi^{\sf nom},\bm{P}^{\sf nom})=100.$ We start by comparing the worst-case performances of $\pi^{\sf nom}$ for $\PP^{\sf (r)}$ and $\PP^{\sf (s)}$.

\begin{table}
\centering
\begin{tabular}{lrrrr}
    \hline
Budget of deviation $\tau$  \;  & \; 0.05 \; & \; 0.07 \; & \; 0.09 \; \\
      \hline
Worst-case of $\pi^{\sf nom}$ for $\PP^{\sf (r)}$  & 94.40 & 92.21 & 90.04  \\
Worst-case of $\pi^{\sf nom}$ for $\PP^{\sf (s)}$ & 91.74 & 88.56 & 85.46 \\
      \hline
    \end{tabular}
          \caption{Comparison of nominal and worst-cases over $\PP^{\sf (r)}$ and $\PP^{\sf (s)}$ for the optimal nominal policy $\pi^{\sf nom}$.}
          \label{tab:PE_nom_K}
\end{table}
The reward of the optimal policy may deteriorate; for instance, for $\tau =0.07$, the worst-case of $\pi^{\sf nom}$ for $\PP^{\sf (r)}$ is $92.21$ and $88.56$ for $\PP^{\sf (s)}$, to compare with $100$ for the nominal kernel $\bm{P}^{\sf nom}.$ Moreover, the set $\PP^{\sf (r)}$ seems to yield a less conservative estimation of the worst-case of $\pi^{\sf nom}$ than the set $\PP^{\sf (s)}$. Indeed, in this example the worst-case of $\pi^{\sf nom}$ for $\PP^{\sf (r)}$ is always higher than the worst-case of $\pi^{\sf nom}$ for $\PP^{\sf (s)}$. Note that $\PP^{\sf (r)}$ is not a subset of $\PP^{\sf (s)},$ because in  $\PP^{\sf (s)}$ the deviation in $\| \cdot \|_{1}$ is constrained to be smaller than $\sqrt{S \cdot A}\cdot \tau,$ whereas there is no such constraint in $\PP^{\sf (r)}$.

We now compute an optimal robust policy $\pi^{\sf rob,r}$ for $\PP^{\sf (r)}$ using Algorithm \ref{alg:VI-PI} and an optimal robust policy $\pi^{\sf rob,s}$ for $\PP^{\sf (s)}$~\citep{Kuhn}. We compare their worst-case performances (for their respective uncertainty sets) and their performances for the nominal transition kernel $\bm{P}^{\sf nom}.$
\begin{table}
\centering
\begin{tabular}{lrrr}
    \hline
Budget of deviation $\tau$   & \; 0.05 \; & \; 0.07 \; & \; 0.09 \\
    \hline
Nominal reward of $\pi^{\sf rob,r}$ & 100.00 & 100.00 & 100.00   \\
Worst-case of $\pi^{\sf rob,r}$ for $\PP^{\sf (r)}$  & 94.40 & 92.21 & 90.04   \\
            \hline
Nominal reward of $\pi^{\sf rob,s}$  & 99.28 & 98.53 & 97.81 \;  \\
Worst-case of $\pi^{\sf rob,s}$ for $\PP^{\sf (s)}$ & 91.90 & 89.09 & 86.62  \\
      \hline
    \end{tabular}
           \caption{Worst-case and nominal performances of the robust policies $\pi^{\sf rob,r}$ and $\pi^{\sf rob,s}$.}
\label{tab:PE_rob_K}
\end{table}
We note that $\pi^{\sf rob, r}$ is identical to $\pi^{\sf nom}$ for $\tau \in \{0.05,0.07,0.09\}$. This indicates that $r$-rectangularity captures sensitivity when necessary, whereas $\pi^{\sf rob,s}$ deviates from $\pi^{\sf nom}$, for only a moderate improvement in worst-case: for $\tau=0.09$, compare the worst-case $86.62$ of $\pi^{\sf rob,s}$ with $85.46$, the worst-case of $\pi^{\sf nom}.$ We also note that in all our experiments the policy $\pi^{\sf rob,s}$ was randomized, which can be hard to interpret and implement in practice.

Since nature might not be adversarial, we compare the performances of $\pi^{\sf rob,r}$ and $ \pi^{\sf rob,s}$ on a sample of kernels around the nominal transitions $\bm{P}^{\sf nom}.$
The robust policy $\pi^{\sf rob,r}$ maximizes the worst-case reward over (some) rank $r$ deviations from the nominal transition kernels $\bm{P}^{\sf nom}$. Therefore, we simulate a random perturbation of rank $r$ from the kernel $\bm{P}^{\sf nom}$ by uniformly generating a random factor matrix and some random coefficients matrices, such that the maximum deviation on each component of the transition kernel is smaller than $\tau$. We also want to consider the case where the coefficients of the perturbations are all independent. 
More precisely, we consider $B_{\sf r}$ and $B_{\infty}$ as
\begin{align*}
B_{\sf r}& = \{ \bm{P} \; | \; \bm{P}_{s}=\bm{P}^{\sf nom}_{s}+\bm{Wu}_{s}, \| \bm{P} - \bm{P}^{\sf nom} \|_{\infty} \leq \tau, \bm{W} \in \R^{S \times r}_{+},(\bm{u}_{s})_{s \in \X} \in \R^{(r \times A) \times S}_{+}, \bm{P}_{sa}^{\top}\bm{e}=1, \forall \; (s,a) \in \X \times \A \}, \\
B_{\infty} &= \{ \bm{P}  \; | \; \| \bm{P} - \bm{P}^{\sf nom} \|_{\infty} \leq \tau, \bm{P}_{sa}^{\top}\bm{e}=1, \forall \; (s,a) \in \X \times \A \}.
\end{align*}
Note that $B_{\infty}$ contains the uncertainty sets $\PP^{\sf (r)}$ and $\PP^{\sf (s)}$. We draw $10000$ kernels $\bm{P}$ uniformly in $B_{\sf r}$ and $B_{\infty}$ and we present in Table \ref{tab:PE_emp_K} the means and $95 \%$ confidence levels $\sf conf_{95}$ of the rewards $R(\pi^{\sf rob,r},\bm{P})$ and $R(\pi^{\sf rob,s},\bm{P})$ for different values of the parameter $\tau >0.$  We would would like to recall that the policy $\pi^{\sf rob,s}$ changes with the parameter $\tau$. 
\begin{table}
    \centering
	\begin{tabular}{lrrrrrr}
    	\hline
    	Budget of deviation $\tau$ & \multicolumn{2}{c}{0.05} & \multicolumn{2}{c}{0.07} & \multicolumn{2}{c}{0.09}  \\
 & \; mean \;  & \; $\sf conf_{95}$ & \; mean \;  & \; $\sf conf_{95}$ & \; mean \;  & \; $\sf conf_{95}$ \\
	      \hline
	     	Empirical reward of $\pi^{\sf rob,r}$ in $B_{\sf r}$ & 99.923 & 0.002 & 99.893 & 0.002 & 99.864 & 0.002   \\
	Empirical reward of $\pi^{\sf rob,s}$ in $B_{\sf r}$ & 99.203 & 0.002 & 98.425 & 0.002 & 97.685 & 0.002  \\
	\hline
		     	Empirical reward of $\pi^{\sf rob,r}$ in $B_{\infty}$ & 98.463 & 0.013 & 97.976 & 0.017 & 97.554 & 0.021  \\
	Empirical reward of $\pi^{\sf rob,s}$ in $B_{\infty}$ & 97.734 & 0.012 & 96.638 & 0.015 & 95.793 & 0.018  \\           
   	   \hline
    \end{tabular}
       \caption{Empirical performances of the policies $\pi^{\sf rob,r}$ and $\pi^{\sf rob,s}$. We draw $10000$ kernels $\bm{P}$ in $B_{\sf r}$ and $B_{\infty}$ and we report the means of the ratio $ \dfrac{R(\pi,\bm{P})}{R(\pi^{\sf nom},\bm{P}^{\sf nom}})$ and the $95 \%$ confidence levels, defined as $1.96 \cdot \text{std}/\sqrt{10000}$ where `std' stands for the standard deviations of the observed rewards.}
       \label{tab:PE_emp_K}
\end{table}

We see empirically that $\pi^{\sf rob,s}$ performs worse than $\pi^{\sf rob,r}$, both in $B_{\sf r}$ and $B_{\infty}.$ For instance, for a maximum deviation of $\tau = 0.09$, the empirical mean of the rewards of $\pi^{\sf rob,r}$ in $B_{\sf r}$ is $99.864$, to compare to $97.685$ for $\pi^{\sf rob,s}.$ In $B_{\infty}$, the empirical mean for $\pi^{\sf rob,r}$ is $97.544$, higher than the mean $95.793$ for $\pi^{\sf rob,s}.$  Moreover, for a same budget of deviation $\tau$, we notice that the mean of the rewards is higher when we sample kernels in $\PP^{\sf (r)}$ than in $\PP^{\sf (s)}$.
\subsection{An example inspired by healthcare.}\label{subsec:second-ex}
 MDPs have been used in healthcare applications, and are efficient to analyze chronic diseases and treatments that may evolve over time, see \cite{mdp-med-1} and \cite{mdp-med-2} for short surveys of the applications of MDP to medical decision making.  
Here, we consider an example inspired from a healthcare application, where we model the evolution of the patient's health using Markovian transition. The decision maker (doctor) prescribes a drug dosage at every state. Transitions from different health states are likely to be related as they are influenced by the same underlying factors, such as genetics, demographics, and/or physiologic characteristics of certain diseases. Therefore, $(s,a)$ and $s$-rectangular may be too conservative in this setting. Note that in a real healthcare example, the MDP would require a large number of states and actions to accurately depict the health dynamics of a patient; this section provides an illustrative example to demonstrate the usefulness of our model compared to other uncertainty sets.

We use $\X = \{1,2,3,4,5,m\}$ where $m$ is an absorbing \textit{mortality} state. We use $\A = \{a_{1} = \text{ low}, a_{2} = \text{ medium}, a_{3} = \text{ high}\}$ for the drug dosage. This goal is to minimize the mortality rate of the patients, while reducing the invasiveness of the treatment prescribed. 
 We consider an estimated nominal transition kernel $\bm{P}^{\sf nom}$ given to us and an upper bound $\tau>0$ on the maximum possible deviation from any component of the nominal kernel. We give the details of the rewards and transition probabilities in Figure \ref{fig:Ex_BoU}.
\begin{figure}[h]
\begin{subfigure}{0.32\textwidth}
  \includegraphics[width=0.9\linewidth]{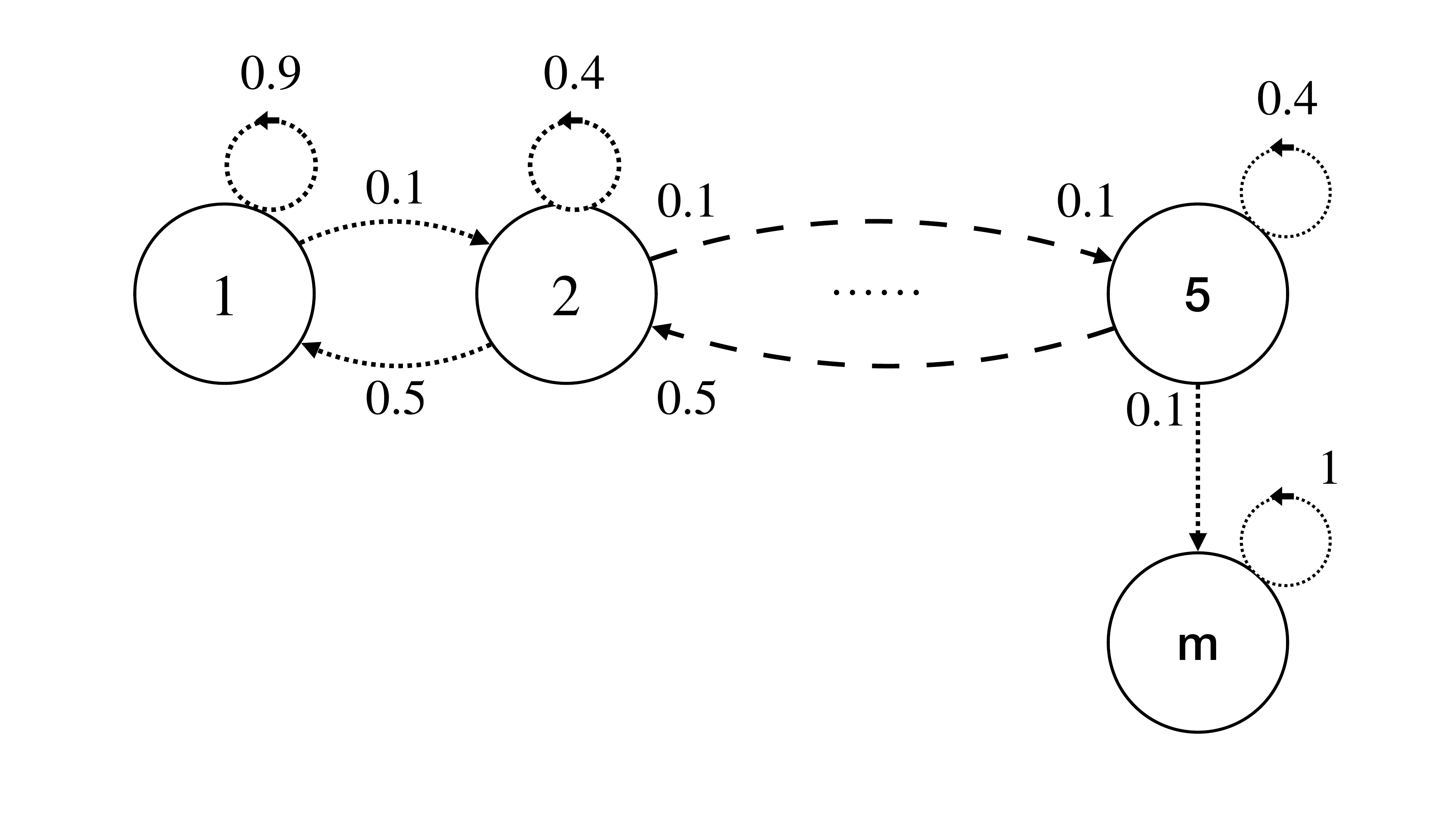}
\caption{Transition probabilities for \\ action {\it high dosage}.}
\label{fig:healthcare-high}
\end{subfigure}
\begin{subfigure}{0.32\textwidth}
  \includegraphics[width=0.90\linewidth]{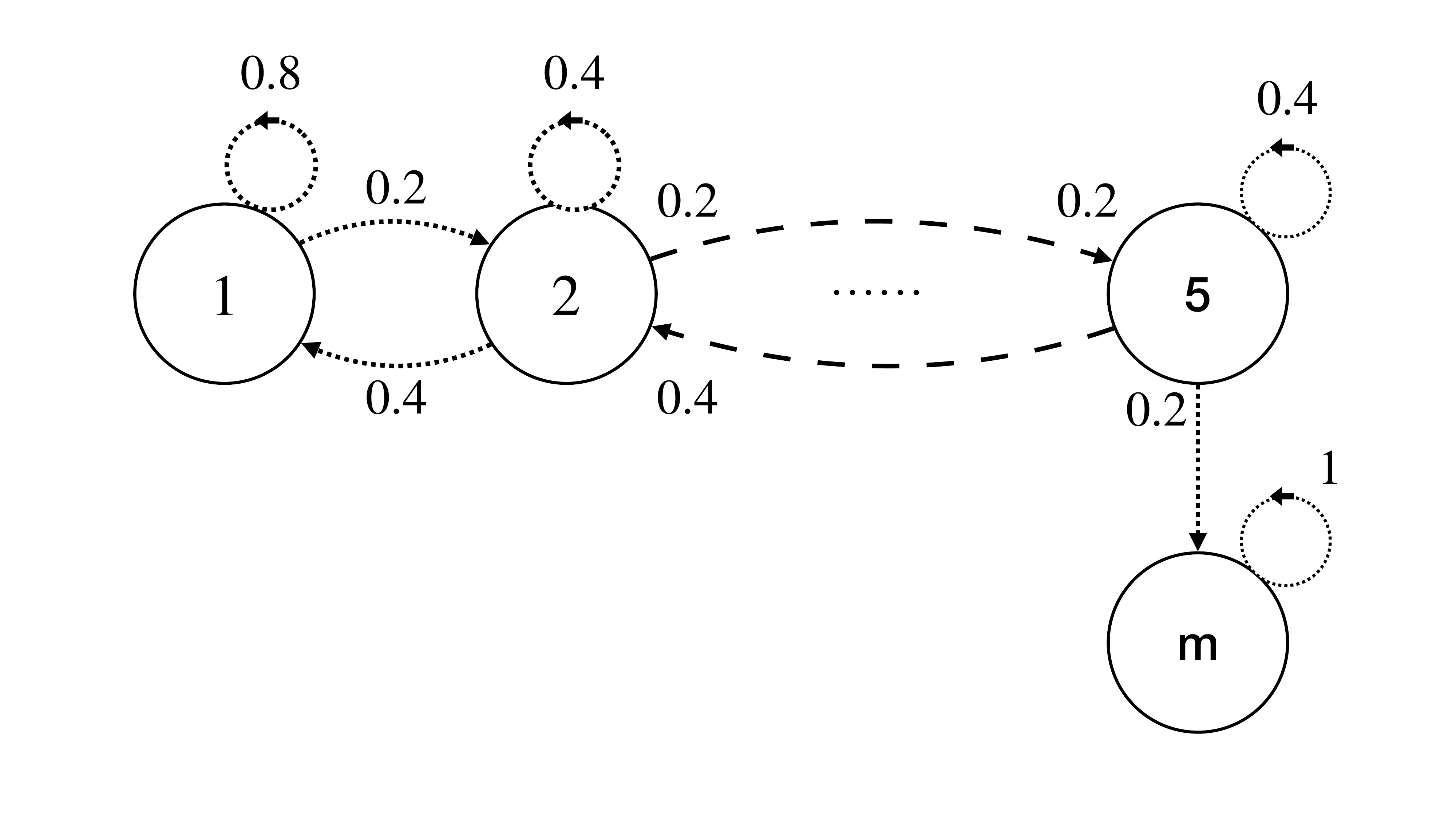}
\caption{Transition probabilities for \\ action {\it medium dosage}.}
\label{fig:healthcare-low}
\end{subfigure}\begin{subfigure}{0.32\textwidth}
  \includegraphics[width=0.9\linewidth]{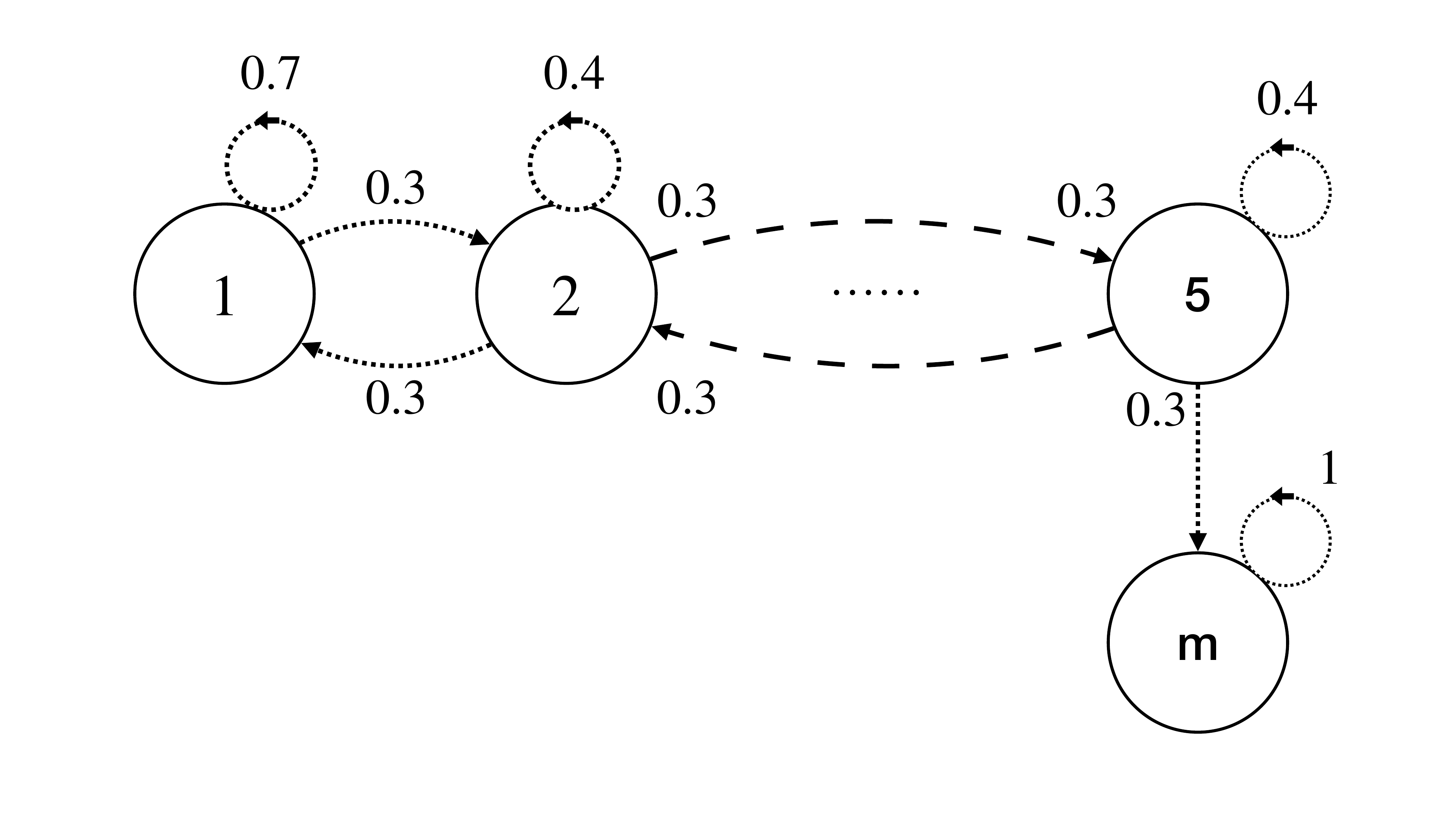}
\caption{Transition probabilities for \\ action {\it low dosage}.}
\label{fig:healthcare-wait}
\end{subfigure}
\caption{Transition probabilities for the healthcare MDP instance. }
\label{fig:Ex_BoU}
\end{figure}
As in the previous example we construct two uncertainty sets $\PP^{\sf (r)}$ and $\PP^{\sf (s)}$ and we compute the worst-case of the optimal nominal policy $\pi^{\sf nom}$. We then compare the performances of the robust policies associated with these two models with the performances of $\pi^{\sf nom}.$

\noindent \textbf{Construction of the $r$-rectangular uncertainty set.} For the $r$-rectangular uncertainty set, we use the construction of the previous example with budget of uncertainty sets.  {\color{black} For this healthcare example, we choose $r=8$ to obtain a good NMF approximation of the nominal transition kernel. This corresponds to a rank reduction of $10$ compared to $S \cdot A = 6 \cdot 3 =18.$}We would like to emphasize that the NMF decomposition enables us to cope with deterministic transitions. In particular, the state $m$ (representing mortality) is absorbing and is not subject to any uncertainty. We model this by increasing the parameter $r$ to $r+1$, introducing $\W^{r+1} = \{ (0,0,0,0,0,1) \}$ and defining $(u^{i}_{ma})_{a \in \A,i=1,...,r+1}$ as the matrix with zero everywhere and $u^{r+1}_{ma} = 1,$ for all action $a \in \A$.

{\color{black}
\noindent \textbf{Construction of the $s$-rectangular uncertainty set.}
We consider the same budget of uncertainty as in the previous example:
\begin{equation*}
\begin{aligned}
\PP^{\sf (s)}_{s} & =\{\bm{P}_{s}=\bm{P}_{s}^{\sf nom} + \bm{\Delta} \; | \;  \bm{\Delta} \in \R^{A \times S}, \; \| \bm{\Delta} \|_{1} \leq \sqrt{S \cdot A} \cdot \tau,\| \bm{\Delta} \|_{\infty} \leq \tau, \; \bm{P}_{s}\bm{e}_{S}=\bm{e}_{A}, \; \bm{P}_{s} \geq \bm{0} \}, s=1,...,5, \\
\PP^{\sf (s)} & = \underset{s \in \X}{\times} \; \PP^{\sf (s)}_{s}.
\end{aligned}
\end{equation*}
The set $\PP^{\sf (s)}_{m}$ reduces to a single matrix since the state $m$ is absorbing:
$$\PP^{\sf (s)}_{m}=\{\bm{P}^{\sf (s)_{m}} \}, \text{ where } P^{\sf (s)}_{m,as'} = 1 \text{ if } s'=m, \text{ and } 0 \text{ otherwise.}$$
}

\noindent \textbf{Empirical results.}
We compute the optimal nominal policy $\pi^{\sf nom}$ using value iteration and we choose the reward parameters such that $R(\pi^{\sf nom},\bm{P}^{\sf nom})=100.$ We compare the worst-case performances of $\pi^{\sf nom}$ for the uncertainty sets $\PP^{\sf (r)}$ and $\PP^{\sf (s)}$ and present our results in Table \ref{tab:PE_rob}.

Note that the performances of $\pi^{\sf nom}$ significantly deteriorates even for small parameter deviations. We observe that the set $\PP^{\sf (r)}$ yields less conservative estimate of the worst-case reward than the set $\PP^{\sf (s)}$. In other words, the worst-case of $\pi^{\sf nom}$ for $\PP^{\sf (r)}$ is always higher than its worst-case for $\PP^{\sf (s)}$.

We also compare the performance of the optimal robust policies $\pi^{\sf rob,r}$ and $\pi^{\sf rob,s}$. The comparisons are presented in Table \ref{tab:PE_rob}. The two robust policies have comparable performances on the nominal kernel $\bm{P}^{\sf nom}.$ The increase in worst-case is proportionally higher for $\pi^{\sf rob,s}$ than for $\pi^{\sf rob,r}$ : for $\tau = 0.09$, the worst-case increase from $31.51$ for $\pi^{\sf nom}$ over $\PP^{\sf (s)}$ to $38.69$ for $\pi^{\sf rob,s}.$ For the same maximum deviation $\tau$, in $\PP^{\sf (r)}$ the worst-case increases from $35.63$ for $\pi^{\sf nom}$ to $36.56$ for $\pi^{\sf rob,r}.$ Yet, this may point out that the policy $\pi^{\sf rob,s}$ sacrifices performance on other feasible kernels in order to increase the worst-case performance.
\begin{table}
\centering
\begin{tabular}{lrrr}
    \hline
Budget of deviation $\tau$   & \;  $\tau=$0.05 \; & \; $\tau=$0.07 \; & \; $\tau=$0.09 \\
    \hline
Nominal reward of $\pi^{\sf nom}$  & 100.00 & 100.00 & 100.00  \\
Worst-case of $\pi^{\sf nom}$ for $\PP^{\sf (r)}$  & 50.26 & 41.74 & 35.63  \\
Worst-case of $\pi^{\sf nom}$ for $\PP^{\sf (s)}$ & 45.75 & 37.37 & 31.51 \\
            \hline
Nominal reward of $\pi^{\sf rob,r}$  & 100.00 & 92.92 & 92.92  \\
Worst-case of $\pi^{\sf rob,r}$ for $\PP^{\sf (r)}$  & 50.26 & 42.29 & 36.56  \\
            \hline
Nominal reward of $\pi^{\sf rob,s}$  & 91.48 & 91.35 &  89.56 \\
Worst-case of $\pi^{\sf rob,s}$ for $\PP^{\sf (s)}$ & 52.09 & 44.39 & 38.69  \\
      \hline
    \end{tabular}
           \caption{Worst-case and nominal performances of the robust policies $\pi^{\sf rob,r}$ and $\pi^{\sf rob,s}$.}
\label{tab:PE_rob}
\end{table}
Since the perturbations are not necessarily adversarial in practice, we compare the performance of $\pi^{\sf rob,r}$ and $ \pi^{\sf rob,s}$ on randomly sampled kernels around the nominal transitions $\bm{P}^{\sf nom}$. In particular we use the same balls $B_{\sf r}$ and $B_{\infty}$ as in the previous example.
We draw $10000$ kernels $\bm{P}$ uniformly in $B_{\sf r}$ and $B_{\sf \infty}$ and compute the means and $95 \%$ confidence levels $\sf conf_{95}$ of the rewards $R(\pi^{\sf rob,r},\bm{P})$ and $R(\pi^{\sf rob,s},\bm{P})$, for different values of the parameter $\tau >0.$  We would would like to emphasize that the policies $\pi^{\sf rob,s}$ and $\pi^{\sf rob,r}$ changes with the parameter $\tau$. 

\begin{table}
    \centering
	\begin{tabular}{lrrrrrr}
    	\hline
    	Budget of deviation $\tau$ & \multicolumn{2}{c}{$\tau=$0.05} & \multicolumn{2}{c}{$\tau=$0.07} & \multicolumn{2}{c}{$\tau=$0.09}  \\
 & \; mean \;  & \; $\sf conf_{95}$ & \; mean \;  & \; $\sf conf_{95}$ & \; mean \;  & \; $\sf conf_{95}$ \\
	      \hline
Empirical reward of $\pi^{\sf rob,r}$ in $B_{\sf r}$ & 98.66 & 0.03 & 91.60 & 0.03 & 91.26 & 0.04  \\
Empirical reward of $\pi^{\sf rob,s}$ in $B_{\sf r}$ & 90.16 & 0.03 & 89.51 & 0.04 & 87.26 & 0.005  \\
	\hline
Empirical reward of $\pi^{\sf rob,r}$ in $B_{\infty}$ & 82.08 & 0.18 & 76.57 & 0.20 & 73.77 & 0.24  \\
Empirical reward of $\pi^{\sf rob,s}$ in $B_{\infty}$ & 74.52 & 0.11 & 70.11 & 0.13 & 64.94 & 0.15  \\           
   	   \hline
    \end{tabular}
       \caption{Empirical performances of the policies $\pi^{\sf rob,r}$ and $\pi^{\sf rob,s}$. We draw $10000$ kernels $\bm{P}$ in $B_{\sf r}$ and $B_{\infty}$ and we report the means of the reward
$R(\pi,\bm{P})/R(\pi^{\sf nom},\bm{P}^{\sf nom})$ and the $95 \%$ confidence levels, defined as $1.96 \cdot \text{std}/\sqrt{10000}$ where `std' stands for the standard deviations of the observed rewards.}
   \label{tab:PE_emp}
\end{table}

We observe that in all our simulations the policy $\pi^{\sf rob,r}$ has better empirical performance than $\pi^{\sf rob,s}$, for $\tau=0.05,0.07$ or $0.09$ and for both $B_{\sf r}$ and $B_{\infty}$. Moreover, we notice again that the empirical means are higher for the uncertainty set $\PP^{\sf (r)}$ than for the uncertainty set $\PP^{\sf (s)}$. Therefore, this suggests that the $r$-rectangular uncertainty model is a less conservative model for uncertainty than the $s$-rectangular model.

\section{Conclusion.}
As MDP parameters are typically determined based on statistical estimations from real-world data, they are inherently subject to misspecification and estimation errors. Uncertainty in transition parameters may have a negative impact on the performances of the optimal nominal policy in practice, and this has been a limiting factor in implementing MDP-based decisions to real-world problems. We consider \textit{factor matrix} uncertainty sets, which constrain the rank of the deviations from the underlying transition kernel and we highlight their generality. Under the assumption that the factors are unrelated, we prove a strong min-max duality result and we provide an efficient algorithm to compute an optimal robust policy. Our approach allows to model robustness to parameter uncertainty in a more realistic way than in previous state-of-the-art approaches, while remaining tractable.
We also extend some classical properties of MDPs to robust MDPs, most notably the robust maximum principle and Blackwell optimality. From a computational standpoint, we also highlight the construction of factor matrix uncertainty sets using nonnegative matrix factorization. Our numerical experiments suggest that when modeling uncertainty, the decision maker should care about the rank of the potential deviations from the nominal parameters, since empirically low-rank deviations are less conservative than independent perturbations on each component of the nominal transition kernel. 

\bibliographystyle{plainnat}
\bibliography{RMDP-MOR}{}

\appendix
\section{Proof of Theorem \ref{th:PE-tractable}.}\label{app:proof-PE}
The proof of Theorem \ref{th:PE-tractable} uses three lemmas. We start with the following contraction lemma.
\begin{lemma}\label{lem:contraction}[\cite{Nilim}, Lemma 2]
Let $\bm{c} \in \R_{+}^{S}$ and $ f: \R_{+}^{S} \rightarrow \R_{+}^{S}$ be a component-wise non-decreasing contraction. Let $\bm{x}^{*}$ be its unique fixed point. Then $\bm{c}^{\top}\bm{x}^{*} = \max_{\bm{x} \in \R_{+}^{S}, \bm{x} \leq f(\bm{x})} \bm{c}^{\top}\bm{x}.$
\end{lemma}
We also need the following reformulation of the policy evaluation problem.
\begin{lemma}\label{prop:reformu-PE}
Let $\PP$ be a factor matrix uncertainty set. Then the policy evaluation problem can be written as follows.
\begin{equation}\label{eq:reform-PE}
z(\pi) = \bm{p}_{0}^{\top}\bm{r}_{\pi} + \lambda \cdot \min_{ \bm{W} \in \W} \; \bm{p}_{0}^{\top}\bm{T}_{\pi}\left( \bm{I} - \lambda \cdot \bm{W}^{\top}\bm{T}_{\pi}\right)^{-1}\bm{W}^{\top} \bm{r}_{\pi}. 
\end{equation}
\end{lemma}
\proof{Proof of Lemma \ref{prop:reformu-PE}.}
Let $\pi$ be a stationary policy and $\bm{P}$ be a transition kernel in $\PP$. From Lemma 5.6.1 in~\cite{Puterman}, the expected infinite horizon discounted reward can be written 
\[
R(\pi,\bm{P}) = \bm{p}_{0}^{\top}(\bm{I} - \lambda \cdot \bm{L}(\pi,\bm{P})^{\top})^{-1}\bm{r}_{\pi},
\]
where $L(\pi,\bm{P}) \in \R^{S \times S}_{+}$ is the transitions kernel of the Markov chain on $\X$ associated with $\pi$ and $\bm{P}$:
\[ L(\pi,\bm{P})_{ss'}= \sum_{a \in \A} \pi_{sa}P_{sas'},  \forall \; (s,s') \in \X \times \X.\]
From the definition of $\bm{T}_{\pi}$ and $\bm{W}$, we have $
\bm{L}(\pi,\bm{P})  = \bm{T}_{\pi}\bm{W}^{\top}.
$
Hence \begin{align*}
    \left( \bm{I}-\lambda \cdot L(\pi,\bm{P}) \right)^{-1}  =
    \left( \bm{I}-\lambda \cdot \bm{T}_{\pi}\bm{W}^{\top} \right)^{-1} 
     = \sum_{k=0}^{\infty} \lambda^{k} \cdot (\bm{T}_{\pi}\bm{W}^{\top})^{k} 
     =  \bm{I} + \lambda \cdot \bm{T}_{\pi}( \bm{I} - \lambda \cdot \bm{W}^{\top}\bm{T}_{\pi})^{-1}\bm{W}^{\top}.
\end{align*}
Therefore,
$
    z(\pi) = \min_{\bm{P} \in \PP} R(\pi,\bm{P})  =\bm{p}_{0}^{\top}\bm{r}_{\pi} +  \min_{\bm{W} \in \W} \; \lambda \cdot \bm{p}_{0}^{\top}\bm{T}_{\pi}\left( \bm{I} - \lambda \cdot \bm{W}^{\top}\bm{T}_{\pi}\right)^{-1}\bm{W}^{\top} \bm{r}_{\pi}.
$
\hfill \qed
\endproof
Finally, we need the following lemma, which introduces the value function $\bm{\beta} \in \R^{r}$ of the adversary in the adversarial MDP.
\begin{lemma}\label{lem:PE2}
Let $\PP$ be a factor matrix uncertainty set. Then
\begin{align}\label{eq:reform05}
\min_{ \bm{W} \in \W} \;  \bm{p}_{0}^{\top}\bm{T}_{\pi}\left( \bm{I} - \lambda \cdot \bm{W}^{\top}\bm{T}_{\pi}\right)^{-1}\bm{W}^{\top} \bm{r}_{\pi}=\min_{\bm{W} \in \W} \max & \;  \bm{p_{0}}^{\top}\bm{T}_{\pi} \bm{\beta}, \\
& \bm{\beta} \leq \bm{W}^{\top}(\bm{r}_{\pi} + \lambda \cdot \bm{T}_{\pi} \bm{\beta}), \\
& \bm{\beta} \in \R^{r}.
\end{align}
\end{lemma}
\proof{Proof of Lemma \ref{lem:PE2}.}
We define $\bm{\beta} \in \R^{r}$ as a function of $\bm{W} \in \W$:
$\bm{\beta} = (\bm{I} - \lambda \cdot \bm{W}^{\top}\bm{T}_{\pi})^{-1}\bm{W}^{\top}\bm{r}_{\pi}.$
The vector $\bm{\beta}$ is the unique solution of the equation:
\begin{equation}\label{eq:bell-beta}
\bm{\beta} = \bm{W}^{\top}\left(\bm{r}_{\pi} + \lambda \cdot \bm{T}_{\pi}\bm{\beta}\right),
\end{equation}
which can be written component-wise:
\begin{equation}\label{eq:bell-beta-comp-wise}
\beta_{i} = \bm{w}_{i}^{\top}(\bm{r}_{\pi} + \lambda \cdot  \bm{T}_{\pi} \bm{\beta}), \forall \; i \in [r].
\end{equation}
Let us call $LHS$ the value of the optimization program on the left-hand side of \eqref{eq:reform05}. We have
\begin{align}
LHS & =
\min_{\bm{W} \in \W} \; \bm{p}_{0}^{\top} \bm{T}_{\pi}(\bm{I} - \lambda \cdot \bm{W}^{\top}\bm{T}_{\pi})^{-1}\bm{W}^{\top}\bm{r}_{\pi} \\
 & = \label{eq:beta05} \min \{ \bm{p}_{0}^{\top} \bm{T}_{\pi} \bm{\beta} \; | \; \bm{W} \in \W, \bm{\beta} \in \R^{r},\bm{\beta} = (\bm{I} - \lambda \cdot \bm{W}^{\top}\bm{T}_{\pi})^{-1}\bm{W}^{\top}\bm{r}_{\pi} \} \\
& = \label{eq:beta1} \min \{\bm{p}_{0}^{\top} \bm{T}_{\pi} \bm{\beta} \; | \; \bm{W} \in \W,  \bm{\beta} \in \R^{r}, \bm{\beta} = \bm{W}^{\top}(\bm{r}_{\pi} + \lambda \cdot \bm{T}_{\pi}\bm{\beta}) \}  \\
& = \label{eq:beta2} \min_{ \bm{W} \in \W} \max \{ \bm{p}_{0}^{\top} \bm{T}_{\pi} \bm{\beta} \; | \; \bm{\beta} \in \R^{r}, \bm{\beta} \leq \bm{W}^{\top}(\bm{r}_{\pi} + \lambda \cdot \bm{T}_{\pi}\bm{\beta}) \}
\end{align}
where \eqref{eq:beta05} follows the definition of the vector $\bm{\beta}$, Equality \eqref{eq:beta1} follows from \eqref{eq:bell-beta} and \eqref{eq:beta2} follows from Lemma \ref{lem:contraction}.
\hfill \qed
\endproof
We are now ready to prove Theorem \ref{th:PE-tractable}.

\proof{Proof of Theorem \ref{th:PE-tractable}.}
Using the reformulation of Lemma \ref{lem:PE2}, the policy evaluation problem becomes
\begin{align}\label{eq:z-pi-reform}
z(\pi)=\min_{\bm{W} \in \W} \max & \;  \bm{p_{0}}^{\top}(\bm{r}_{\pi} + \lambda \cdot \bm{T}_{\pi} \bm{\beta}), \\
& \bm{\beta} \leq \bm{W}^{\top}(\bm{r}_{\pi} + \lambda \cdot \bm{T}_{\pi}\bm{\beta}), \\
& \bm{\beta} \in \R^{r}.
\end{align}
The gist of the proof of Theorem \ref{th:PE-tractable} is to show that $z(\pi)=\hat{z}(\pi)$, where
\begin{align}\label{eq:RHS}
 \hat{z}(\pi)=\max & \;  \bm{p_{0}}^{\top}(\bm{r}_{\pi} + \lambda \cdot \bm{T}_{\pi} \bm{\beta}) \\
& \; \beta_{i} \leq \min_{\bm{w_{i}} \in \W^{i}}  \bm{w_{i}}^{\top}(\bm{r}_{\pi} + \lambda \cdot \bm{T}_{\pi}\bm{\beta}), \; \forall \; i \in [r] , \label{eq:RHS:tight} \\
& \bm{\beta} \in \R^{r}.
\end{align}
Because of Lemma \ref{lem:contraction}, at optimality in \eqref{eq:RHS} each of the constraint \eqref{eq:RHS:tight} is tight. Let $\bm{\beta}^{a}$ be the solution of \eqref{eq:RHS} and $\bm{W}^{a}=(\bm{w}_{1}^{a},...,\bm{w}_{r}^{a})$ the factor matrix that attains each of the minimum on the components of $\bm{\beta}^{a}$:
\[ \beta^{a}_{i}= \bm{w}_{i}^{a \; \top}(\bm{r}_{\pi} + \lambda \cdot \bm{T}_{\pi}\bm{\beta}^{a}) =\min_{\bm{w}_{i} \in \W^{i}}  \; \bm{w}_{i}^{\top}(\bm{r}_{\pi} + \lambda \cdot \bm{T}_{\pi}\bm{\beta}^{a}) , \; \forall \; i \in [r].\]
These equations uniquely determine the vector $\bm{\beta}^{a}$ since 
\begin{align*}
 \forall \; i \in [r], \beta^{a}_{i}= \bm{w}_{i}^{a \; \top}(\bm{r}_{\pi} + \lambda \cdot \bm{T}_{\pi}\bm{\beta}^{a}) & \iff  \bm{\beta}^{a}= \bm{W}^{a \; \top}(\bm{r}_{\pi} + \lambda \cdot \bm{T}_{\pi} \bm{\beta}^{a}) \\
 & \iff \bm{\beta}^{a} = ( \bm{I} - \lambda \cdot \bm{W}^{a \; \top}\bm{T}_{\pi})^{-1}\bm{W}^{a \; \top}\bm{r}_{\pi}.
\end{align*} 
Note that $\bm{W}^{a}$ is a feasible factor matrix in $\W$ because of the $r$-rectangularity assumption. Therefore, the pair $(\bm{\beta}^{a},\bm{W}^{a})$ is feasible in the optimization problem defining \eqref{eq:z-pi-reform} and $z(\pi) \leq \hat{z}(\pi).$

Following Lemma \ref{lem:contraction}, we also know that the optimum of the program \begin{align*}
 \max & \;  \bm{p_{0}}^{\top}(\bm{r}_{\pi} + \lambda \cdot \bm{T}_{\pi} \bm{\beta}), \\
& \bm{\beta} \leq \bm{W}^{a \; \top}(\bm{r}_{\pi} + \lambda \cdot \bm{T}_{\pi}\bm{\beta}), \\
& \bm{\beta} \in \R^{r},
\end{align*} is attained at a vector $\bm{\hat{\beta}}$ such that $\bm{\hat{\beta}} = \bm{W}^{a \; \top}(\bm{r}_{\pi} + \lambda \cdot \bm{T}_{\pi}\bm{\hat{\beta}}).$ But we just proved that this equation uniquely determines $\bm{\beta}^{a}$ and therefore $\bm{\hat{\beta}}=\bm{\beta}^{a}$.
The matrix $\bm{W}^{a}$ bridges the gap between $z(\pi)$ and $\hat{z}(\pi)$ and the two optimization problems have the same optimum values.

We conclude that
\begin{align}\label{eq:beta4}
 z(\pi) =  \max & \;  \bm{p_{0}}^{\top}(\bm{r}_{\pi} + \lambda \cdot \bm{T}_{\pi} \bm{\beta}) \\
& \; \beta_{i} \leq \min_{\bm{w_{i}} \in \W^{i}}  \bm{w_{i}}^{\top}(\bm{r}_{\pi} + \lambda \cdot \bm{T}_{\pi}\bm{\beta}), \; \forall \; i \in [r], \nonumber \\
& \bm{\beta} \in \R^{r}. \nonumber
\end{align} 
Now let $\phi: \R^{r}_{+} \rightarrow \R^{r}_{+}$ be such that $\phi(\bm{\beta})_{i} = \min_{\bm{w}_{i} \in \W^{i}}  \bm{w}_{i}^{\top}(\bm{r}_{\pi} + \lambda \cdot\bm{T}_{\pi}\bm{\beta}), \forall \; i \in [r].$
Note that $\phi$ is component-wise non-decreasing because of the non-negativity of the sets $\W^{1},...,\W^{r}$ and of the fixed matrix $\bm{T}_{\pi}$. Moreover, $\phi$ is a contraction. Indeed, let $\bm{\beta}_{1}, \bm{\beta}_{2} \in \R^{r}_{+}$ and $i \in [r]$. We have
\begin{align*}
\phi(\bm{\beta}_{1})_{i} & = \min_{\bm{w}_{i} \in \W^{i}}  \bm{w}_{i}^{\top}(\bm{r}_{\pi} + \lambda \cdot\bm{T}_{\pi}\bm{\beta}_{1}) \\ 
& = \min_{\bm{w}_{i} \in \W^{i}}  \bm{w}_{i}^{T}(\bm{r}_{\pi} + \lambda \cdot\bm{T}_{\pi}\bm{\beta}_{2}) + \lambda \cdot \bm{w}_{i}^{\top}\bm{T}_{\pi}(\bm{\beta}_{1}-\bm{\beta}_{2})  \\
& \geq \min_{\bm{w}_{i} \in \W^{i}}  \bm{w}_{i}^{\top}(\bm{r}_{\pi} + \lambda \cdot \bm{T}_{\pi}\bm{\beta}_{2}) + \lambda \cdot \min_{\bm{w}_{i} \in \W^{i}}  \bm{w}_{i}^{\top}\bm{T}_{\pi}(\bm{\beta}_{1}-\bm{\beta}_{2}) \\
& \geq \phi(\bm{\beta_{2}})_{i} + \lambda \cdot \min_{\bm{w}_{i} \in \W^{i}}  \bm{w}_{i}^{\top}\bm{T}_{\pi}(\bm{\beta}_{1}-\bm{\beta}_{2}).
\end{align*}
Therefore, for all $i \in [r]$,
\begin{align*}
\phi(\pi,\bm{\beta}_{2})_{i}-\phi(\pi,\bm{\beta}_{1})_{i}   \leq \lambda \cdot \min_{\bm{w}_{i} \in \W^{i}}  \bm{w}_{i}^{\top}\bm{T}_{\pi}(\bm{\beta}_{2}-\bm{\beta}_{1}) 
 \leq \lambda \cdot \| \bm{\beta}_{2}-\bm{\beta}_{1} \|_{\infty},
\end{align*}
where the last equality follows from $\bm{w}_{i}^{\top}\bm{T}_{\pi}\bm{e}_{r}  =\bm{w}_{i}^{\top}\bm{e}_{S}=1.$
We can do the same computation exchanging the role of $\bm{\beta}_{1}$ and $\bm{\beta}_{2}$, and we conclude that $\phi$ is a contraction.

Therefore, we can apply Lemma \ref{lem:contraction} to the reformulation \eqref{eq:beta4} and we can solve the policy evaluation problem by computing the fixed-point $\bm{\beta}^{*}$ of the contraction $\phi(\pi,\cdot)$, i.e, by computing $\bm{\beta}^{*}$ such that
\begin{equation}\label{Bell-beta-star}
 \beta_{i}^{*} = \min_{\bm{w_{i}} \in \W^{i}}  \bm{w_{i}}^{\top}(\bm{r}_{\pi} + \lambda \cdot \bm{T}_{\pi}\bm{\beta}^{*}), \; \forall \; i \in [r].
\end{equation}
 This can be done by iterating the function $ \phi(\pi,\cdot)$, and Algorithm \ref{alg:VI-PE} is a value iteration algorithm that returns the fixed-point of $\phi(\pi,\cdot)$. From~\cite{Puterman}, the condition $\|\bm{\beta}^{k+1} - \bm{\beta}^{k}\|_{\infty} < \epsilon(1-\lambda)(2 \lambda)^{-1}$ is sufficient to ensure that $ \| \bm{\beta}^{k+1} - \bm{\beta}^{*} \|_{\infty} \leq \epsilon.$ Therefore, Algorithm \ref{alg:VI-PE} returns an $\epsilon$-optimal solution to the policy evaluation problem. We now present the analysis of the running time of Algorithm \ref{alg:VI-PE}.

\vspace{2mm}
\noindent \tb{{\textbf{Running time of Algorithm \ref{alg:VI-PE}.}}}
\begin{itemize}
\item \tb{\textit{Computational cost of initialization.}
Before running Algorithm \ref{alg:VI-PE}, we compute $\bm{r}_{\pi} \in \R^{S}$ and $\bm{T}_{\pi} \in \R^{S \times r}$. This can be done in $O(rSA)$.}
\item \tb{\textit{Number of iterations of Algorithm \ref{alg:VI-PE}.} Following Theorem 6.3.3 in \cite{Puterman}, we can stop Algorithm \ref{alg:VI-PE} as soon as $\|\bm{\beta}^{k+1} - \bm{\beta}^{k}\|_{\infty} < \epsilon(1-\lambda)(2 \lambda)^{-1}$ to ensure that $\| \bm{\beta}^{k+1} - \bm{\beta}^{*} \|_{\infty} \leq \epsilon$, and this condition will be satisfied after at most 
$O ( \log(\epsilon^{-1}))$ iterations of Algorithm \ref{alg:VI-PE}.}
\item \tb{\textit{Computational cost of each iteration.}
 To compute $\bm{\beta}^{k+1}$ given $\bm{\beta}^{k}$,
we compute an $\epsilon_{2}$-solution to the optimization programs $\min_{\bm{w}_{i} \in \W^{i}} \bm{w}_{i}^{\top} \left( \bm{r}_{\pi}+\lambda \bm{T}_{\pi}\bm{\beta}^{k} \right)$, for $i \in [r]$, with $\epsilon_{2}=\epsilon(1-\lambda)(4 \lambda)^{-1}.$ With the notations of Assumption \ref{ass:tractable-u-set}, this can be done in $O\left( \sum_{i=1}^{r} comp\left(\W^{i}\right)\log(\epsilon^{-1}) \right)$ arithmetic operations.
Computing the cost vector $\bm{r}_{\pi} + \lambda \bm{T}_{\pi} \bm{\beta}^{k}$ can be done once at the beginning of each iteration in $O(rS)$.}
\tb{
\noindent
Therefore, the total computational cost for each iteration is in $O\left(rS +  \sum_{i=1}^{r} comp\left(\W^{i}\right)\log(\epsilon^{-1})\right)$. }
\end{itemize}
\tb{
We conclude that Algorithm \ref{alg:VI-PE} stops after a number of \tb{arithmetic operations} of at most $O\left(rSA + rS \log(\epsilon^{-1}) + \sum_{i=1}^{r} comp(\W^{i}) \log^{2}\left(\epsilon^{-1}\right)\right).$}
\hfill \qed
\endproof

\vspace{2mm}
\noindent \textbf{Role of rectangularity of $\W$.}
 As argued in Section \ref{sec:PE}, the Cartesian product property of $\W$ is crucial in our proof of Theorem \ref{th:PE-tractable}. In particular, we provide here an example of an $r$-rectangular uncertainty set, where adding a constraint across the factors $\bm{w}_{1}, ... ,\bm{w}_{r}$ would lead to a different optimal solution.  In particular, consider the following simple MDP instance $\M$:
\begin{itemize}
\item There are two states $s_{1}$ and $s_{2}$ and one action.
\item There is a reward of $0$ in state $s_{1}$ and a reward of $1$ in state $s_{2}$.
\item The number of factors is $r=2$.  The set $\W$ is $\W = \W^{1} \times \W^{2}$ with \[\W^{1} = \Delta(2),\W^{2}= \{ (w_{2,1},w_{2,2}) \in \Delta(2) \; | \; w_{2,1} \leq 1/2\}.\]
\item The coefficients $\left(u_{sa}^{i}\right)_{i,sa}$ are given as
\[  u_{s_{1}}^{1}=1, u_{s_{1}}^{2}=0,u_{s_{2}}^{1}=0,u_{s_{2}}^{2} =1.\]
\item The discount factor can be any $\lambda \in (0,1)$.
\end{itemize}
For this robust MDP instance, the adversarial MDP $\M_{\sf adv}$ is as follows:
\begin{itemize}
\item The states $1$ and $2$ correspond to the value of the index of the factors $\bm{w}_{1}, \bm{w}_{2}$.
\item For $i \in \{1,2\}$,  in state $i$ the adversary picks $\bm{w}_{i} \in \W^{i}$.  
\item The reward $\hat{r}_{i,s}$ for $i \in \{1,2\}$ and $s \in \{s_{1},s_{2}\}$ are simply $\hat{r}_{i,s} = \delta_{\{s=s_{2}\}}$.
\item In terms of transitions, note that we have constructed this instance such that the probability of transitioning from $i \in \{1,2\}$ to $j \in \{1,2\}$ is $w_{i,j}$.
\item The discount factor is the same as in the nominal MDP $\M$.
\end{itemize}
It is straightforward that in this adversarial MDP instance, the optimal policy for the adversary is to choose $\bm{w}_{1}^{*} = (1,0),\bm{w}_{2}^{*} = (1/2,1/2)$ (recall that the adversarial MDP is a \textit{minimization} problem). 
If we add the constraint $\bm{w}_{1}=\bm{w}_{2}$ to $\M_{\sf adv}$ , the previous solution $\left( \bm{w}_{1}^{*},\bm{w}^{*}_{2} \right)$ is not feasible anymore. This highlights the importance of $r$-rectangularity and of the Cartesian product structure of the set $\W$.
Note that for the new robust MDP with an additional constraint, the Bellman equation \eqref{Bell-beta-star} may not hold in the first place, since the new set of factors is not a Cartesian product anymore. In particular, we only proved that \eqref{Bell-beta-star} holds when the set of possible factors $\W$ is a Cartesian product. 
\section{Proof of Lemma \ref{lem:PEfinal}.}\label{app:PE:reform-general}

This follows directly from combining \eqref{eq:beta4} with Lemma \ref{lem:contraction}.
\section{Proof of Proposition \ref{prop:stat-mrkv}.}\label{app:stat-mrkv}
Remember that $\Pi^{G}$ denotes the set of all policies (possibly history-dependent and non-stationary) and $\Pi \subset \Pi^{G}$ is the set of stationary policies.
We want to prove:
$ \max_{ \pi \in \Pi^{G}} \; \min_{\bm{P} \in \PP } \; R(\pi,\bm{P}) = \max_{ \pi \in \Pi} \; \min_{\bm{P} \in \PP } \; R(\pi,\bm{P}).$
We proved in Theorem \ref{th:duality} that
$ \max_{ \pi \in \Pi} \; \min_{\bm{P} \in \PP } \; R(\pi,\bm{P}) =  \min_{\bm{P} \in \PP } \; \max_{ \pi \in \Pi} \; R(\pi,\bm{P}).$
Therefore,
\begin{align}
\max_{ \pi \in \Pi} \; \min_{\bm{P} \in \PP } \; R(\pi,\bm{P}) & \leq \max_{ \pi \in \Pi^{G}} \; \min_{\bm{P} \in \PP } \; R(\pi,\bm{P}) \label{eq:mrkv1} \\
& \leq \min_{\bm{P} \in \PP } \;  \max_{ \pi \in \Pi^{G}} \; R(\pi,\bm{P}) \label{eq:mrkv2} \\
& = \min_{\bm{P} \in \PP } \;  \max_{ \pi \in \Pi} \; R(\pi,\bm{P} )\label{eq:mrkv3}   \\
& =   \max_{ \pi \in \Pi} \; \min_{\bm{P} \in \PP } \; R(\pi,\bm{P}) \label{eq:mrkv4}  
\end{align}
where \eqref{eq:mrkv1} follows from $\Pi \subset \Pi^{G}$, \eqref{eq:mrkv2} follows from weak duality, \eqref{eq:mrkv3} follows from the fact that for a non-robust Markov decision process, an optimal policy can always be found in the set of stationary policies, and \eqref{eq:mrkv4} follows from Theorem \ref{th:duality}.
Therefore all these inequalities are equalities, and
$ \max_{ \pi \in \Pi^{G}} \; \min_{\bm{P} \in \PP } \; R(\pi,\bm{P}) = \max_{ \pi \in \Pi} \; \min_{\bm{P} \in \PP } \; R(\pi,\bm{P}).$

\section{A counter-example to the equilibrium result of Theorem \ref{th:duality}.}\label{app:counter-example-equilibrium}
\tb{
In this appendix, we show that there exists an $s$-rectangular MDP instance where the optimal solutions of the left-hand side and the right-hand side of \eqref{eq:strong-min-max-duality} are not attained by the same pairs.}

\tb{
We consider the same example as in Proposition 1 in \cite{Kuhn}. The MDP is presented in Figure \ref{fig:counter_example}. We
consider an infinite horizon MDP with three states $\{ s_{1},s_{2},s_{3}\}$ and two actions $\{ a_{1},a_{2}\}$. The discount factor is any $\lambda \in (0,1)$. Let $\xi \in [0,1]$.
\begin{itemize}
\item The decision maker always starts in state $s_{1}$.
\item If action $a_{1}$ is selected in state $s_{1}$, the decision maker transitions to state $s_{2}$ with probability $\xi$, and to state $s_{3}$ with probability $1-\xi$. 
\item If action $a_{2}$ is selected in state $s_{1}$, the decision maker transitions to state $s_{3}$ with probability $1-\xi$, and to state $s_{3}$ with probability $\xi$. 
\item The transitions outside of $s_{2}$ and $s_{3}$ do not depend on $a_{1}$ and $a_{2}$, and the states $s_{2}$ and $s_{3}$ are absorbing.
\item There is no reward associated with state $s_{1}$:
\[ r_{s_{1}a_{1}} = r_{s_{2}a_{2}}=0.\]
The rewards for states $s_{2}$ and $s_{3}$ are as follows:
\begin{align*}
r_{s_{2}a_{1}} & = r_{s_{2}a_{2}} = 1, \\
r_{s_{3}a_{1}} & = r_{s_{3}a_{2}} = 0.
\end{align*}
\end{itemize}}
\begin{figure}[h]
\begin{subfigure}{0.5\textwidth}
  \includegraphics[width=0.7\linewidth]{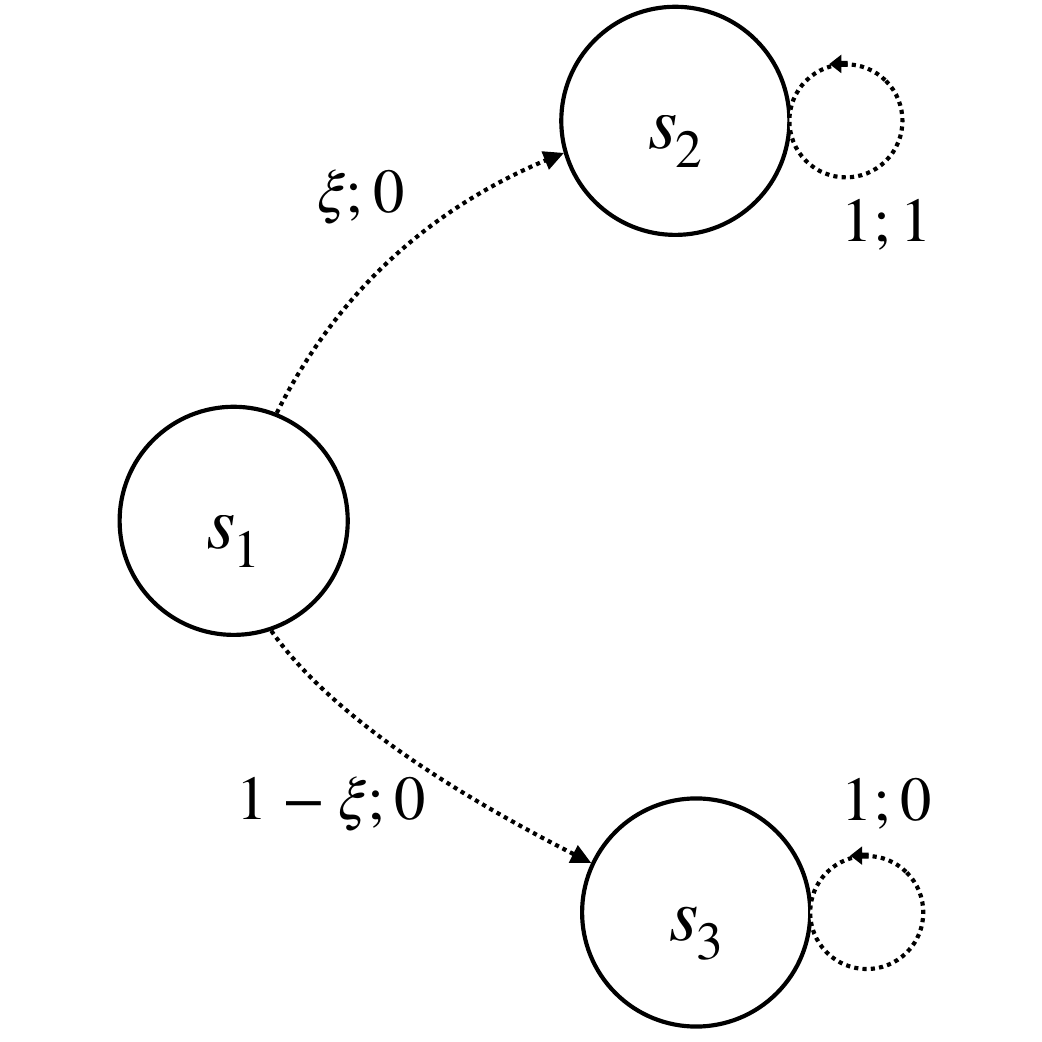}
\caption{Transitions and rewards for choosing action $a_{1}$.}
\label{fig:counter_a_1}
\end{subfigure}
\begin{subfigure}{0.5\textwidth}
  \includegraphics[width=0.7\linewidth]{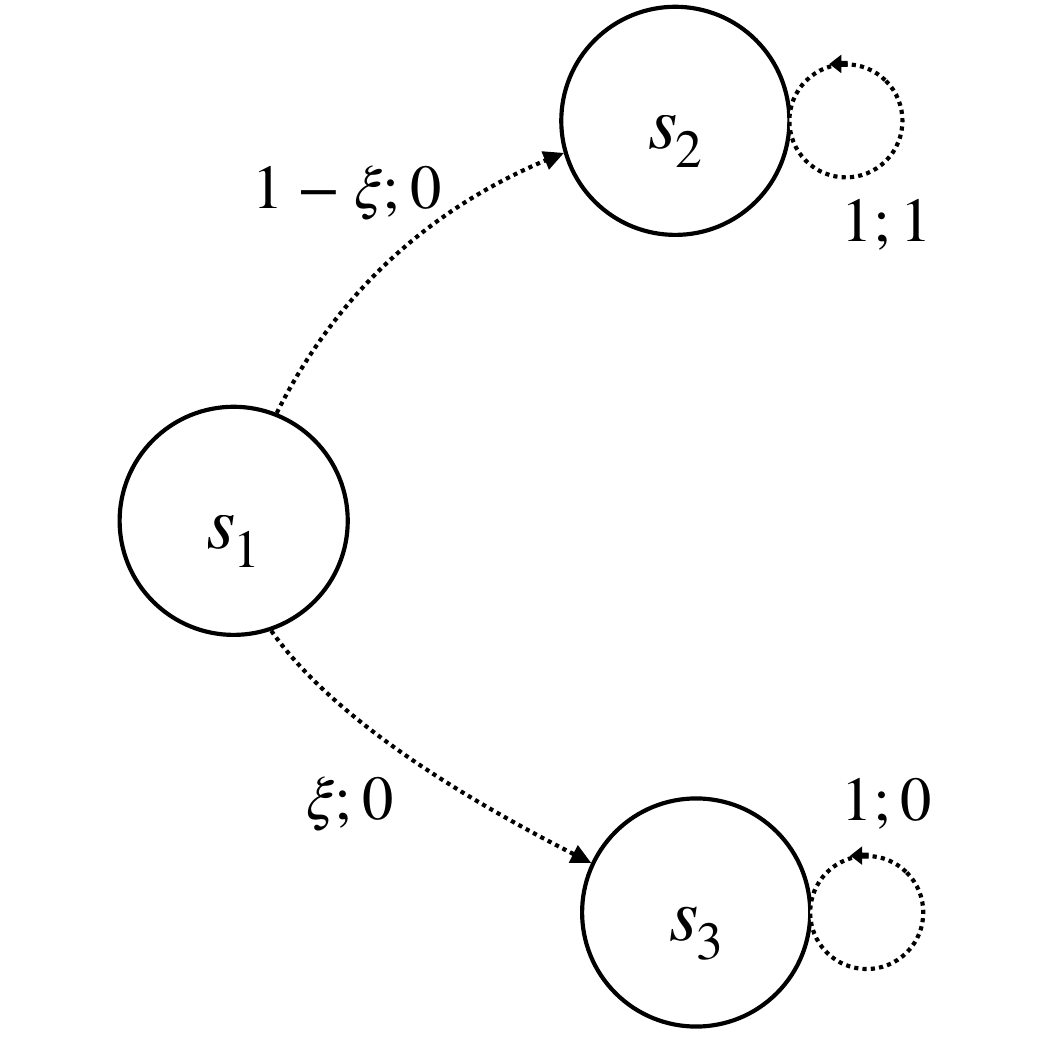}
\caption{Transitions and rewards for choosing action $a_{2}$.}
\label{fig:counter_a_2}
\end{subfigure}
\caption{Counter-example for the equilibrium result of Theorem \ref{th:duality} for an $s$-rectangular MDP. Transition rates and associated rewards are noted above the transition arcs. The expressions before the semicolon correspond to transition probabilities toward the next state, the expressions after the semicolon correspond to the associated rewards.}
\label{fig:counter_example}
\end{figure}

\tb{
For this example, the reward can be parametrized as a function of $\alpha \in [0,1]$ (for $\alpha = \pi_{s_{1}a_{1}}$) and of $\xi \in [0,1]$. We have 
\[ R(\alpha,\xi) = \frac{1}{1-\lambda} \left( \alpha \xi + (1-\alpha)(1-\xi) \right).\]
Consider the $s$-rectangular uncertainty set associated with $\xi \in [0,1/3]$.}

\tb{
On the one hand, we have that
\begin{align*}
\max_{\alpha \in [0,1]} \min_{\xi \in [0,1/3]} \; R(\alpha,\xi) & = \max_{\alpha \in [0,1]} \min_{\xi \in [0,1/3]} \; \frac{1}{1-\lambda} \left( \alpha \xi + (1-\alpha)(1-\xi) \right) \\
& = \frac{1}{1-\lambda} \max_{\alpha \in [0,1]} \min \{\frac{2}{3}-\frac{\alpha}{3},1-\alpha \} \\
& =  \frac{1}{1-\lambda} \frac{2}{3},
\end{align*}
attained for $\alpha_{\sf left}=0,\xi_{\sf left} = 1/3$.}

\tb{
On the other hand, we have that 
\begin{align*}
 \min_{\xi \in [0,1/3]} \max_{\alpha \in [0,1]} \; R(\alpha,\xi) & =  \min_{\xi \in [0,1/3]} \max_{\alpha \in [0,1]} \; \frac{1}{1-\lambda} \left( \alpha \xi + (1-\alpha)(1-\xi) \right) \\
& = \frac{1}{1-\lambda} \min_{\xi \in [0,1/3]} \max \{\xi,1-\xi \} \\
& = \frac{1}{1-\lambda} \frac{2}{3},
\end{align*}
attained for $\alpha_{\sf right}=1,\xi_{\sf right} = 1/3$.}

\tb{
Therefore, for this particular $s$-rectangular MDP instance, different policies are solving the $\max \min$ and the $\min \max$ optimization problems, even though strong duality holds.}
%
%
%

\section{Proof of Lemma \ref{lem:F1}.}\label{app:lem:F1}
 Let $F:\R_{+}^{S} \rightarrow \R_{+}^{S}$ be such that
\[F(\bm{v})_{s} = \max_{\pi_{s} \in \Delta } \{ r_{\pi_{s},s} +  \lambda \cdot (\bm{T}_{\pi}(\min_{\bm{w}_{i} \in \W^{i}} \{ \bm{w}_{i}^{\top}\bm{v} \})_{i \in [r]})_{s} \}, \forall \; s \in \X, \forall \; \bm{v} \in \R_{+}^{S}.\]
We prove here that this function is a component-wise non-decreasing contraction. 
\begin{lemma}\label{lem:F1} The function $F$ is a component-wise non-decreasing contraction.
\end{lemma}
\proof{Proof.}
The mapping  $F$ is component-wise non-decreasing because $\W \subseteq \R^{S \times r}$ and $\bm{T}_{\pi} \in \R^{S \times  r}$ have non-negative entries.

Let us prove that $F$ is a contraction.
Let $\bm{v}_{1}$ and $\bm{v}_{2}$ be in $\R^{S}_{+}$ and $s \in \X$.
\begin{align*}
F(\bm{v}_{1})_{s} & = \max_{\pi_{s}\in \Delta } \{ r_{\pi_{s},s} +  \lambda \cdot (\bm{T}_{\pi}(\min_{\bm{w}_{i} \in \W^{i}} \{ \bm{w}_{i}^{\top}\bm{v}_{1} \})_{i \in [r]})_{s} \} \\
& = \max_{\pi_{s} \in \Delta } \{ r_{\pi_{s},s} +  \lambda \cdot (\bm{T}_{\pi}(\min_{\bm{w}_{i} \in \W^{i}} \{ \bm{w}_{i}^{\top}\bm{v}_{2}+\bm{w}_{i}^{\top}(\bm{v}_{1}-\bm{v}_{2}) \})_{i \in [r]})_{s} \} \\
& \leq \max_{\pi_{s}\in \Delta } \{ r_{\pi_{s},s} +  \lambda \cdot (\bm{T}_{\pi}(\min_{\bm{w}_{i} \in \W^{i}} \{ \bm{w}_{i}^{\top}\bm{v}_{2}+ \|\bm{v}_{1}-\bm{v}_{2}\|_{\infty}\})_{i \in [r]})_{s} \} \\
& \leq \max_{\pi_{s} \in \Delta} \{ r_{\pi_{s},s} +  \lambda \cdot (\bm{T}_{\pi}(\min_{\bm{w}_{i} \in \W^{i}} \{ \bm{w}_{i}^{\top}\bm{v}_{2}\})_{i \in [r]})_{s} \} + \lambda \cdot \|\bm{v}_{1}-\bm{v}_{2}\|_{\infty} \\
& \leq F(\bm{v}_{2})_{s} + \lambda \cdot \|\bm{v}_{1}-\bm{v}_{2}\|_{\infty}.
\end{align*}
Therefore, for all $s \in \X$, $F(\bm{v}_{1})_{s} - F(\bm{v}_{2})_{s} \leq \lambda \cdot \|\bm{v}_{1}-\bm{v}_{2}\|_{\infty}.$
Inverting the role of the two vectors we can conclude that $F$ is a contraction:
$\| F(\bm{v}_{1}) - F(\bm{v}_{2})\|_{\infty} \leq \lambda \cdot \|\bm{v}_{1}-\bm{v}_{2}\|_{\infty}.$
\hfill \qed
\section{Proof of Proposition \ref{prop:robust-max-principle}.}\label{app:robust-max-principle}
We call $\bm{v}^{\pi,\bm{W}}$ the value function of the decision maker associated with the policy $\pi$ and the factor matrix $\bm{W}$.
\begin{enumerate}
\item  Let $\pi$ be a policy and $\bm{W}^{1} \in \arg \min_{\bm{W} \in \W} \; R(\pi,\bm{W}).$
Let $\bm{W}^{0} \in \W$ be any factor matrix. From the Bellman Equation \eqref{eq:bell-beta} and Lemma \ref{lem:Wv-beta}, we have
$\bm{v}^{\pi,\bm{W}^{0}}=\bm{r}_{\pi} + \lambda \cdot \bm{T}_{\pi}\bm{\beta}^{0},\text{ where }
\bm{\beta}^{0}=\bm{W}^{0 \; \top}(\bm{r}_{\pi} + \lambda \cdot \bm{T}_{\pi}\bm{\beta}^{0}).$ From the Bellman Equation \eqref{Bell-beta-star} and Lemma \ref{lem:Wv-beta},
$$\bm{v}^{\pi,\bm{W}^{1}}=\bm{r}_{\pi} + \lambda \cdot \bm{T}_{\pi}\bm{\beta}^{1}, \text{ where }\bm{\beta}^{1}=\bm{W}^{1 \; \top}(\bm{r}_{\pi} + \lambda \cdot \bm{T}_{\pi}\bm{\beta}^{1}) = \left(\min_{\bm{w}_{i} \in \W^{i}} \bm{w}_{i}^{\top}(\bm{r}_{\pi} + \lambda \cdot \bm{T}_{\pi}\bm{\beta}^{1})\right)_{i \in [r]}.$$
From Theorem \ref{th:PE-tractable}, we know that the sequence of vectors $(\phi^{n}(\pi,\bm{\beta}^{0}))_{n \in \N}$ converges to $\bm{\beta}^{1}$. Moreover, for any $n \in \N$, we have the component-wise inequality:
$\bm{\beta}^{1} \leq \phi^{n}(\pi,\bm{\beta}^{0}) \leq \bm{\beta}^{0}.$
The matrix $\bm{T}_{\pi}$ being non-negative, we obtain $\bm{r}_{\pi} + \lambda \cdot \bm{T}_{\pi}\bm{\beta}^{1} \leq \bm{r}_{\pi} + \lambda \cdot \bm{T}_{\pi}\bm{\beta}^{0} ,$
and we conclude that
$\bm{v}^{\pi,\bm{W}^{1}} \leq \bm{v}^{\pi,\bm{W}^{0}}. $

\item Let $\pi$ be a policy and let $\bm{W}^{1} \in \arg \min_{\bm{W} \in \W} \; R(\pi,\bm{W}).$ We write $\left( \bm{T}_{\pi_{s}}\left(\min_{\bm{w}_{i} \in \W^{i}} \bm{w}_{i}^{\top}\bm{v}^{\pi,\bm{W}^{1}}\right)_{i \in [r]} \right)_{s}$ for the $s$-th component of the vector $ \bm{T}_{\pi}\left(\min_{\bm{w}_{i}\in \W^{i}} \bm{w}_{i}^{\top}\bm{v}^{\pi,\bm{W}^{1}}\right)_{i \in [r]} \in \R^{S}_{+}.$

The Bellman equation for the value function $\bm{v}^{\pi,\bm{W}^{1}}$ yields, for each $s \in \X$,
\begin{align*} v^{\pi,\bm{W}^{1}}_{s}  = r_{\pi_{s},s} + \lambda \cdot \left( \bm{T}_{\pi_{s}}\left(\min_{\bm{w}_{i}\in \W^{i}} \bm{w}_{i}^{\top}\bm{v}^{\pi,\bm{W}^{1}}\right)_{i \in [r]} \right)_{s} 
&  \leq r_{\pi_{s},s} + \lambda \cdot \left( \bm{T}_{\pi_{s}}\left( \bm{w}_{i}^{* \; \top}\bm{v}^{\pi,\bm{W}^{1}}\right)_{i \in [r]} \right)_{s}  \\
 & \leq \max_{\bm{\tilde{\pi}}_{s} \in \Delta} r_{\tilde{\pi}_{s},s} + \lambda \cdot \left( \bm{T}_{\tilde{\pi}_{s}}\left( \bm{w}_{i}^{* \; \top}\bm{v}^{\pi,\bm{W}^{1}}\right)_{i \in [r]} \right)_{s}
,\end{align*}
which can be written component-wise as:
\begin{equation}\label{eq:r-max-princi}
\bm{v}^{\pi,\bm{W}^{1}} \leq H_{\bm{W}^{*}}(\bm{v}^{\pi,\bm{W}^{1}}),
\end{equation}
where $H_{\bm{W}^{*}}(\cdot): \R^{S}_{+} \rightarrow \R^{S}_{+}$ is defined as $H_{\bm{W}^{*}}(\bm{v})_{s} = \max_{\bm{\tilde{\pi}}_{s} \in \Delta} r_{\tilde{\pi}_{s},s} + \lambda \cdot \left( \bm{T}_{\tilde{\pi}_{s}}\left( \bm{w}_{i}^{* \; \top}\bm{v}\right)_{i \in [r]} \right)_{s}.$
 $H_{\bm{W}^{*}}$ is a non-decreasing contraction and its fixed-point is the value function of the optimal policy for the MDP with the transition kernel associated with the factor matrix $\bm{W}^{*}$. From Theorem \ref{th:duality},
$$ (\pi^{*},\bm{W}^{*}) \in \arg \max_{\pi \in \Pi} \min_{\bm{W} \in \W} R(\pi,\bm{W}) \iff (\bm{W}^{*},\pi^{*}) \in \arg \min_{\bm{W} \in \W} \max_{\pi \in \Pi} R(\pi,\bm{W}) .$$
Therefore, $\pi^{*}$ is the optimal policy for the MDP with the transition kernel associated with $\bm{W}^{*}$.
Iterating \eqref{eq:r-max-princi}, using the fact that $H_{\bm{W}^{*}}$ is non-decreasing and that for any initial vector $\bm{v}_{0}$, $\lim_{n \rightarrow \infty} H_{\bm{W}^{*}}^{n}(\bm{v}_{0}) = \bm{v}^{\pi^{*},\bm{W}^{*}},$
we can conclude that for all policy $\pi$, 
$\bm{W}^{1} \in \arg \; \min_{\bm{W} \in \W} R(\pi,\bm{W}) \Rightarrow v^{\pi,\bm{W}^{1}}_{s} \leq v^{\pi^{*},\bm{W}^{*}}_{s}, \forall \; s \in \X.$
\end{enumerate}
\section{Proof of Blackwell optimality.}
\tb{
We start with the following lemma, that will be useful for proving Proposition \ref{prop:black} and Proposition \ref{prop:blackwell-01}. We provide a concise proof of Lemma \ref{lem:tiroir-infinity} for the sake of completeness.
\begin{lemma}\label{lem:tiroir-infinity}
Let $(u_{n})_{n \geq 0}$ be a sequence with values in a finite set $\EE$. Then there exists an element $e \in \EE$ that is attained infinitely often by $(u_{n})_{n \geq 0}$.
\end{lemma}
\proof{Proof.}
For each element $e \in \EE$, let $\NN_{e}=
\{ n \in \N \; | \; u_{n}=e \}$. Note that $\N = \bigcup_{e \in \EE} \NN_{e}$, because $\forall \; n \in \N, u_{n} \in \EE$. Therefore, if the sets $\NN_{e}$ were all finite, then $\N$ would be finite, as the finite union of some finite sets. This implies that there exists at least one element $e \in \EE$ for which $\NN_{e}$ is infinite.
\hfill \qed
\endproof
}

\subsection{Proof of Proposition \ref{prop:black}.}\label{app:black}
\tb{We start by providing a high-level overview of our proof. We proceed by contradiction and assume that  there does not exist a discount factor $\hat{\lambda} \in (0,1)$ for which the {\it same} pair of policy-kernel is optimal for all $\lambda \in [\hat{\lambda},1)$. In this case, the optimal pair of policy-kernel is varying infinitely often as the discount factor $\lambda$ approaches $1$. This is will result in a contradiction with the fact that both the set of policies and kernels have finitely many extreme points (under the assumption of Proposition \ref{prop:black}), and that the value function $\lambda \mapsto \bm{v}_{\lambda}^{\pi,\bm{P}}$ is a continuous {\it rational} function on $(0,1)$, for a fixed policy $\pi$ and a fixed factor matrix $\bm{W}$.}

\tb{
We now prove Proposition \ref{prop:black}.}
We call $\bm{v}^{\pi,\bm{W}}_{\lambda}$ the value function of the decision maker associated with the policy $\pi$, the factor matrix $\bm{W}$ and the discount factor $\lambda$. We call $R(\pi,\bm{W},\lambda)$ the reward \eqref{eq-expreward} associated with the policy $\pi$, the factor matrix $\bm{W}$ and the discount factor $\lambda.$
We call $z^{*}_{\lambda}$ the policy improvement problem with discount factor $\lambda$:
$z^{*}_{\lambda} = \max_{\pi \in \Pi} \; \min_{\bm{W} \in \W} \; R(\pi,\bm{W},\lambda).$
Let $(\pi_{\lambda},\bm{W}_{\lambda})$ be a solution of $z^{*}_{\lambda}$ with $\pi_{\lambda}$ stationary and deterministic. From \eqref{eq:intric-fix-pt}, the factor matrix $\bm{W}_{\lambda}$ belongs to the set of extreme points of $\W$. Since the set of stationary deterministic policies and the set of extreme points of $\W$  are finite, \tb{from Lemma \ref{lem:tiroir-infinity}} we can choose $(\lambda_{n})_{n \geq 0}$ such that there exists a fixed $(\pi^{*},\bm{W}^{*})$ such that
$$ \lambda_{n} \rightarrow 1, \text{ and } (\pi^{*},\bm{W}^{*}) \in \arg \max_{\pi \in \Pi} \; \min_{\bm{W} \in \W} \; R(\pi,\bm{W},\lambda_{n}),  \forall \; n \geq 0 .$$
We prove that the pair $(\pi^{*},\bm{W}^{*})$ is an optimal solution to $z^{*}_{\lambda}$ for all discount factor $\lambda$ sufficiently close to $1$.

Let us assume the contrary, i.e., let us assume that there exists a sequence of discount factor $(\gamma_{n})_{n \geq 0}$ such that $$\gamma_{n} \rightarrow 1, \text{ and }(\pi^{*},\bm{W}^{*}) \notin \arg \max_{\pi \in \Pi} \; \min_{\bm{W} \in \W} \; R(\pi,\bm{W},\gamma_{n}), \forall \; n \geq 0.$$ From the finiteness of the set of stationary deterministic policies and of the set of extreme points of $\W$, \tb{from Lemma \ref{lem:tiroir-infinity}} we can choose $(\tilde{\pi},\bm{\tilde{W}})$ such that this pair is optimal in $z^{*}_{\gamma_{n}}$ for all $n \geq 0.$

Now for all $n \in \N$,  $(\pi^{*},\bm{W}^{*}) \notin \arg \max_{\pi \in \Pi} \; \min_{\bm{W} \in \W} \; R(\pi,\bm{W},\gamma_{n}) \Rightarrow z_{\gamma_{n}}(\pi^{*}) < z_{\gamma_{n}}(\tilde{\pi}),$
and the robust maximum principle of Proposition \ref{prop:robust-max-principle} implies that for all $n \in \N$, there exists a state $x_{1,n}$ and a factor matrix $\bm{W}^{*,n}$ such that the value functions satisfy
$v_{\gamma_{n},x_{1,n}}^{\pi^{*},\bm{W}^{*,n}} < v_{\gamma_{n},x_{1,n}}^{\tilde{\pi},\bm{\tilde{W}}}.$
From the finiteness of the set of extreme points of $\W$ and of the set of states $\X,$ \tb{from Lemma \ref{lem:tiroir-infinity}} we can chose $x_{1}$ and $\bm{W}^{**}$ such that for any $n \in \N$,
$v_{\gamma_{n},x_{1}}^{\pi^{*},\bm{W}^{**}} < v_{\gamma_{n},x_{1}}^{\tilde{\pi},\bm{\tilde{W}}}.$
Similarly, the robust maximum principle gives, for any $n \in \N$, $$\pi^{*} \in \arg \max_{\pi \in \Pi} \; \min_{\bm{W} \in \W} \; R(\pi,\bm{W},\lambda_{n}) \Rightarrow  z_{\lambda_{n}}(\tilde{\pi}) \leq z_{\lambda_{n}}(\pi^{*}) \Rightarrow \exists \; \bm{\tilde{W}}^{n} \in \W, v_{\lambda_{n},x_{1}}^{\tilde{\pi},\bm{\tilde{W}}^{n}} \leq v_{\lambda_{n},x_{1}}^{\pi^{*},\bm{W}^{*}}. $$
From the finiteness of the set of extreme points of $\W$, \tb{from Lemma \ref{lem:tiroir-infinity}} we can choose all the factor matrices $(\bm{\tilde{W}}^{n})_{n \in \N}$ to be equal to the same factor matrix $\bm{\tilde{\tilde{W}}}$. 
Therefore, there exists a state $x_{1}$ such that for any integer $n$,
$v_{\gamma_{n},x_{1}}^{\pi^{*},\bm{W}^{**}} < v_{\gamma_{n},x_{1}}^{\tilde{\pi},\bm{\tilde{W}}},$
 $ v_{\lambda_{n},x_{1}}^{\tilde{\pi},\bm{\tilde{\tilde{W}}}} \leq v_{\lambda_{n},x_{1}}^{\pi^{*},\bm{W}^{*}}. $
 Now from the robust maximum principle,
 $\bm{W}^{*} \in  \arg \min_{\bm{W} \in \W} R(\pi^{*},\bm{W},\lambda_{n}) \Rightarrow v_{\lambda_{n},x_{1}}^{\pi^{*},\bm{W}^{*}} \leq v_{\lambda_{n},x_{1}}^{\pi^{*},\bm{W}^{**}},$
 and for the same reason
 $v_{\gamma_{n},x_{1}}^{\tilde{\pi},\bm{\tilde{W}}} \leq v_{\gamma_{n},x_{1}}^{\tilde{\pi},\bm{\tilde{\tilde{W}}}}.$
Therefore,
 $v_{\lambda_{n},x_{1}}^{\tilde{\pi},\bm{\tilde{\tilde{W}}}} \leq v_{\lambda_{n},x_{1}}^{\pi^{*},\bm{W}^{**}}.$
Overall, we have constructed two factor matrices $\bm{W}^{**}$ and $\bm{\tilde{\tilde{W}}}$ such that
\begin{equation}\label{eq:contradiction}
 v_{\gamma_{n},x_{1}}^{\pi^{*},\bm{W}^{**}} < v_{\gamma_{n},x_{1}}^{\tilde{\pi},\bm{\tilde{\tilde{W}}}}, v_{\lambda_{n},x_{1}}^{\tilde{\pi},\bm{\tilde{\tilde{W}}}}  \leq v_{\lambda_{n},x_{1}}^{\pi^{*},\bm{W}^{**}}.
\end{equation}
Let us define the function $f:(0,1) \rightarrow \R, f(t) = v_{t,x_{1}}^{\tilde{\pi},\bm{\tilde{\tilde{W}}}} - v_{t,x_{1}}^{\pi^{*},\bm{W}^{**}}$.
For any stationary policy $\pi$, any kernel $\bm{P}$, any discount factor $t \in (0,1),$ the value function $\bm{v}^{\pi,\bm{P}}_{t}$ satisfies
 $\bm{v}^{\pi,\bm{P}}_{t}= (\bm{I}-t \bm{L}(\pi,\bm{P}))^{-1}\bm{r}_{\pi}, \text{ with } L(\pi,\bm{P})_{ss'} = \sum_{a \in \A} \pi_{sa}P_{sas'}, \forall \; (s,s') \in \X \times \X.$ Therefore Cramer's rule implies that the function $f$ is a continuous rational function on $(0,1)$, i.e., it is the ratio of two polynomial of finite degrees and the denominator does not have any zeros in $(0,1)$.

But \eqref{eq:contradiction} implies that the function $f$ takes the value $0$ for an infinite number of scalars $\theta_{n} \rightarrow 1$, and non-zero values for an infinite number of scalars $\gamma_{n} \rightarrow 1.$
A continuous rational function cannot take the value 0 at an infinite number of different points $\{ \theta_{n} \; | \; n \geq 0 \}$ without being itself the zero function. This is a contradiction, and therefore the pair $(\pi^{*},\bm{W}^{*})$ is an optimal solution of $z^{*}_{\lambda}$ for all $\lambda$ sufficiently close to $1.$
{\color{black}
\subsection{Proof of Proposition \ref{prop:blackwell-01}.}\label{app:black-01}
The proof of Proposition \ref{prop:black} extends to any $\lambda \in (0,1]$. In particular, for any $\lambda \in (0,1]$, there exists an open interval $\I_{\lambda}$ around $\lambda$ such that the same pair of policy-factor matrix is optimal on $\I_{\lambda}$ in the robust MDP problem.  The open intervals $\{ \I_{\lambda} \; | \; \lambda \in (0,1] \}$ form a covering of the compact set $[0,1]$. From the Borel-Lebesgue covering theorem, there exists a \textit{finite} numbers of $\lambda_{0}^{'}, ..., \lambda_{m}^{'}$ such that $\bigcup_{i=0}^{m} \I_{\lambda_{i}^{'}} = [0,1]$, for some $m \in \N$. From this, for some $p \in \N$, we can construct $\lambda_{0}=0 < \lambda_{1} < ... < \lambda_{p}=1$ such that the same pair of stationary deterministic policy and factor matrix $(\pi_{i},\bm{W}_{i})$ is an optimal solution to the robust MDP problem on each interval $(\lambda_{j},\lambda_{j+1})$ for all $j \in \{0,...,p-1\}$.}

\end{document}